\documentclass[final,a4paper,fleqn,10pt]{cas-sc}

\usepackage[numbers,sort&compress]{natbib}

\usepackage[debug]{Table}
\usepackage[below]{placeins}
\usepackage{afterpage}

\usepackage{paralist}
\usepackage{empheq}

\DeclareMathOperator{\cond}{cond}
\DeclareMathOperator{\supp}{supp}
\newcommand{\du}{\ensuremath{\,\mathrm{d}}}
\newcommand{\ttimes}{$\times$}

\newcommand{\coltext}{(For interpretation of the references to colour in this figure legend, the reader is referred to the web version of this article.)}
\usepackage{hyperref}
\hypersetup{
    pdftitle =    {Imposing nonlocal boundary conditions in Galerkin-type methods based on non-interpolatory functions},
    pdfauthor =  {Svajunas Sajavicius; Thomas Takacs},
    pdfsubject =  {Computers and Mathematics with Applications, 80 (2020) 2877-2895. doi:10.1016/j.camwa.2020.09.016},
    pdfkeywords = {Nonlocal boundary conditions; Galerkin methods; Non-interpolatory basis functions; Isogeometric analysis},
    pdfcreator =  {},
    pdfproducer = {},
    linkcolor = {hscolor},
    urlcolor =  {hscolor},
    citecolor = {hscolor},
    pdfstartview =  FitH,
    bookmarksopen = True,
    colorlinks =    True
}
\usepackage[open, openlevel=3]{bookmark}

\makeatletter
\long\def\@makecaption#1#2{%
   \vskip\abovecaptionskip
   \sbox\@tempboxa{\footnotesize\textbf{#1.} #2}%
   \ifdim \wd\@tempboxa >\hsize
     \footnotesize{\textbf{#1.} #2}\par
   \else
     \global \@minipagefalse
     \hb@xt@\hsize{\hfil\box\@tempboxa\hfil}%
   \fi
   \vskip\belowcaptionskip}
\makeatother

\begin{document}
\let\WriteBookmarks\relax
\def\floatpagepagefraction{1}
\def\textpagefraction{.001}
\shorttitle{S. Sajavi\v{c}ius and T.~Takacs \hfill Computers and Mathematics with Applications 80 (2020) 2877--2895}
\shortauthors{\href{https://doi.org/10.1016/j.camwa.2020.09.016}{https://doi.org/10.1016/j.camwa.2020.09.016}}

\title [mode = title]{Imposing nonlocal boundary conditions in {G}alerkin-type methods based on
  non-interpolatory functions}                      

\author[1]{Svaj\={u}nas Sajavi\v{c}ius}
\cormark[1]
\ead{svajunas.sajavicius@ktu.lt}
\ead[url]{https://www.personalas.ktu.lt/~svasaj/}
\address[1]{
  Department of Software Engineering, Faculty of Informatics,
  Kaunas University of Technology, Lithuania
}
\credit{Conceptualization, Methodology, Software, Investigation, Writing - original draft}

\author[2]{Thomas Takacs}
\ead{thomas.takacs@jku.at}
\ead[URL]{http://www.ag.jku.at/~thomast/}
\address[2]{
  Institute of Applied Geometry,
  Johannes Kepler University Linz, Austria
}
\credit{Validation, Formal analysis, Writing - review \& editing}

\cortext[cor1]{Corresponding author}

\begin{abstract}
The imposition of inhomogeneous Dirichlet (essential) boundary conditions 
 is a fundamental challenge in the application of Galerkin-type methods based 
 on non-interpolatory functions, i.e., functions which do not possess the 
 Kronecker delta property. Such functions typically are used in various 
 meshfree methods, as well as methods based on the isogeometric paradigm. The 
 present paper analyses a model problem consisting of the Poisson equation 
 subject to non-standard boundary conditions. Namely, instead of classical 
 boundary conditions, the model problem involves Dirichlet- and Neumann-type 
 nonlocal boundary conditions. Variational formulations with strongly and 
 weakly imposed inhomogeneous Dirichlet-type nonlocal conditions are derived 
 and compared within an extensive numerical study in the isogeometric 
 framework based on non-uniform rational B-splines (NURBS). The attention in 
 the numerical study is paid mainly to the influence of the nonlocal 
 boundary conditions on the properties of the considered discretisation 
 methods.
\end{abstract}


\ExplSyntaxOn
\keys_set:nn { stm / mktitle } { nologo }
\ExplSyntaxOff
\maketitle

\section{Introduction}

The problem of imposing inhomogeneous Dirichlet (essential) boundary 
conditions is well-known in the context of Galerkin-type meshfree methods. 
Such methods, as the element-free Galerkin (EFG) 
method~\cite{Belytschko1994}, the reproducing kernel particle method 
(RKPM)~\cite{Liu1995} or the corrected smooth particle hydrodynamics (CSPH) 
method~\cite{Bonet1999} use basis functions which do not satisfy the 
Kronecker delta property, i.e.~a basis function associated with a particle 
does not necessarily vanish at other particles. Due to the non-interpolatory 
nature of the basis functions, the imposition of inhomogeneous Dirichlet 
boundary conditions in the meshfree approaches is not as straightforward as, 
for example, in the finite element method based on the classical 
interpolatory functions.

There exist many techniques for the imposition of classical essential 
boundary conditions in meshfree methods. For example, the applications of 
Lagrange multipliers~\cite{Belytschko1994}, modified variational 
principles~\cite{Belytschko1994,Dumont2006}, penalty 
methods~\cite{Zhu1998,Bonet2000,Cho2008}, perturbed 
Lagrangian~\cite{Chu1995}, modifying meshfree basis 
functions~\cite{Gosz1996,Guenther1998,Chen2000,Wagner2000,Sukumar2004,Oh2009}, 
or coupling meshfree approaches and mesh-based techniques (namely, finite 
elements)~\cite{Belytschko1995,Krongauz1996,Huerta2000,Wagner2001} can be 
mentioned. A novel method to treat Neumann and Dirichlet boundary conditions 
in meshfree methods for elliptic equations using nodal integration is 
presented in~\cite{Fougeron2019}.  

Non-interpolatory basis functions, such as non-uniform rational B-splines 
(NURBS), are also used in isogeometric analysis 
(IGA)~\cite{Hughes2005,Cottrell2009}. IGA is a computational technique which 
has been proposed with the aim to establish a direct connection between 
computer-aided design (CAD) and engineering analysis. Currently, the methods 
based on the isogeometric paradigm receive a lot of attention in the 
literature. A review of the mathematical foundation and recent results 
related to IGA can be found 
in~\cite{BeiraodaVeiga2014,Hughes2017,Hughes2018}. Some trends and 
developments of IGA, as well as computer implementation aspects, are 
summarised in~\cite{Nguyen2015}.

The imposition of inhomogeneous essential boundary conditions in the IGA 
framework is considered frequently (see 
e.g.~\cite{Costantini2010,Wang2010b,Chen2011b,Mitchell2011,Govindjee2012}). A 
strong (direct) imposition is a possible way to treat such conditions, but 
weak imposition methods often are more efficient. For example, a weak 
imposition of essential boundary conditions is superior to strong imposition 
when solving various problems in fluid 
mechanics~\cite{Bazilevs2007,Bazilevs2007b}. In several situations considered 
in isogeometric analysis the imposition of boundary or coupling conditions is 
non-trivial, e.g.,~when the domain is composed of trimmed 
patches~\cite{Ruess2014} or when the boundary of the physical domain is 
described as a curve or surface immersed in a regular 
grid~\cite{Schillinger2012,Rank2012,Nitti2020}. 

Weak methods for the imposition of inhomogeneous essential boundary 
conditions in the context of meshfree methods usually are based on a certain 
modification of the basis functions (see 
e.g.~\cite{Belytschko1995,Gosz1996,Huerta2000}) or on a modification of the 
variational formulation of the problem. The well-known examples of methods 
from the latter class are the Lagrange multiplier method, the penalty method 
and Nitsche's method. The Lagrange multiplier method is general and can be 
applied to the solution of various problems. However, it has some significant 
disadvantages, cf.~\cite{Fernandez-Mendez2004,Huerta2017}, among them is the 
difficulty of choosing a suitable Lagrange multiplier space. The penalty and 
Nitsche's methods are attractive alternatives to the Lagrange multiplier 
method. A Nitsche formulation for spline-based finite elements has been 
applied to second- and fourth-order problems in~\cite{Embar2010}. 

The present paper deals with a model problem for Poisson equation with 
non-standard boundary conditions. Instead of the classical Dirichlet or 
Neumann boundary conditions, we consider a model problem involving 
\emph{nonlocal boundary conditions} which can be used in order to express 
various nonlocal couplings between different parts of the model geometry. 
Such conditions, sometimes also called \emph{side conditions}, appear in 
various mathematical models. For example, the models arising in 
thermoelasticity~\cite{Day1985b}, thermodynamics~\cite{Day1992}, biological 
fluid dynamics~\cite{Hazanee2014}, plasma physics~\cite{Diaz1998}, 
geology~\cite{Glotov2016} or peridynamics~\cite{Madenci2018} should be 
mentioned. Nonlocal multipoint constraints are used in the finite element 
analysis related to structural mechanics~\cite{Ainsworth2001,Jendele2009}. A 
steady Poisson equation endowed with a generalised Robin boundary condition 
involving a Laplace--Beltrami operator on the boundary has been analysed 
in~\cite{Kashiwabara2015}. Overlapping Schwarz methods are based on a 
coupling the boundary of one subdomain with the interior of another 
subdomain, which is then solved iteratively. This approach has been applied 
in the isogeometric context in~\cite{Bercovier2015}. In~\cite{Kargaran2019}, 
this coupled problem on overlapping subdomains (called overlapping 
multi-patch structure) is solved all at once. Hence, the coupling, which is 
performed point-wise at Greville abscissae, can be interpreted as an 
imposition of nonlocal constraints on parts of the boundary of each subdomain. 

The main focus of this paper is the imposition of nonlocal boundary 
conditions for Galerkin methods using non-interpolatory basis functions. We 
investigate computationally the effects that such conditions have on the 
properties of the methods. While the numerical study is performed in the IGA 
framework, at least some of the insights can be applied in the context of 
other discretisation approaches too. 

The paper is organised as follows. In Section~\ref{sec:2}, we formulate the model 
problem with nonlocal boundary conditions. Several variational formulations 
of the model problem are given in  Section~\ref{sec:3}. We apply methods for the 
imposition of inhomogeneous Dirichlet-type nonlocal boundary conditions in 
the IGA context. Therefore, in  Section~\ref{sec:4} we shortly review the basic 
concepts and ideas used in IGA. The discretised equations are given in 
Section~\ref{sec:5}, and the results of the numerical study are presented in 
Section~\ref{sec:6}. Finally, we conclude the paper in Section~\ref{sec:7} with a 
summary and several remarks. 

\section{Model problem}

\label{sec:2}

Partial differential equations (PDEs) with nonlocal boundary conditions 
receive a lot of attention in the literature from the point of view of both 
their theoretical investigation (see 
e.g.~\cite{Avalishvili2011,Avalishvili2018}) and their numerical analysis 
(see~\cite{Sajavicius2014a,Sajavicius2016} and references therein). 

In this paper we consider, as a model problem, the two-dimensional Poisson 
equation with mixed nonlocal boundary conditions: Find a function $u: 
\overline{\Omega} \rightarrow \mathbb{R}$, such that
\begin{equation}
\label{eq:1}
\left\{\begin{array}{lll} -\nabla \cdot (\kappa \nabla u) = f & \text{in $\Omega$,} \rule{8cm}{0pt} & \mbox{(a)}\\
  u + \mathcal{L}_1[u] = g & \text{on $\Gamma_D$,}  & \mbox{(b)}\\
  \kappa \nabla u \cdot \mathbf{n} + \mathcal{L}_2[u] = h & \text{on $\Gamma_N$,}  & \mbox{(c)} \end{array}\right.
\end{equation}
where $\Omega \subset \mathbb{R}^2$ is an open and bounded domain with 
Lipschitz boundary $\partial \Omega \equiv \Gamma = \overline{\Gamma_D \cup 
\Gamma_N}$, $\Gamma_D \cap \Gamma_N = \emptyset$, $\Gamma_D \neq \emptyset$, 
$\mathbf{n}$ is the outward normal unit vector on $\partial \Omega$, $\kappa: 
\Omega \rightarrow \mathbb{R}$, $f: \Omega \rightarrow \mathbb{R}$, $g: 
\Gamma_D \rightarrow \mathbb{R}$ and $h: \Gamma_N \rightarrow \mathbb{R}$ are 
given functions, $\mathcal{L}_1: H^1(\Omega) \rightarrow \mathbb{R}$ and 
$\mathcal{L}_2: H^1(\Omega) \rightarrow \mathbb{R}$ are given linear 
functionals. We refer to (\hyperref[eq:1]{1b}) and (\hyperref[eq:1]{1c}) as \emph{Dirichlet-} 
and \emph{Neumann-type nonlocal boundary conditions}, respectively. 

The linear functionals $\mathcal{L}_1$ and $\mathcal{L}_2$ define nonlocal 
parts of the boundary conditions. Various functionals can be used in order to 
express different types of nonlocal couplings. For example, the functional
\begin{equation}
\label{eq:2}
    \mathcal{L}[u] = \sum_{l}^{}{\gamma_l u(\mathbf{x}_l^*)}
\end{equation}
for $\mathbf{x}_l^* \in \Omega$ represents \emph{discrete coupling}, while
the functional
\begin{equation}
\label{eq:3}
    \mathcal{L}[u] = \int_{\Omega}{\gamma u}\,\du\Omega
\end{equation}
corresponds to more general \emph{integral coupling}. The weights $\gamma_l$ 
and $\gamma$ are constants or functions measuring (controlling) the influence 
of the functionals.

Even though point evaluations are not well-defined in $H^1$, that means they 
are not suitable for the variational formulation we introduce in the 
following section, we can nonetheless mimic a point evaluation at 
$\mathbf{x}^*_l$ by introducing a suitable integral coupling. Let 
$\rho^{\delta}_{l} \in C^\infty(\Omega)$ be such that $\supp 
\rho^{\delta}_{l} \subseteq V_{\delta}(\mathbf{x}^*_l) \cap \Omega$, where 
$V_{\delta}(\mathbf{x}^*_l)$ is the $\delta$-neighbourhood of 
$\mathbf{x}^*_l$, and
\[
    \int_{\Omega}{\rho^{\delta}_{l}}\du\Omega = \gamma_l.
\]
Then the integral coupling functional as in~\eqref{eq:3}
\[
    \mathcal{L}^{\delta}[u] = \sum_l \int_{\Omega}{\rho^{\delta}_{l} u}\du\Omega\]
is well-defined for $u \in H^1(\Omega)$ and approximates the discrete 
coupling functional in~\eqref{eq:2}. If $u$ is sufficiently smooth, i.e.~$u\in H^{1+\varepsilon}(\Omega)$, then $\lim_{\delta \rightarrow 0} 
\mathcal{L}^{\delta}[u]$ is also well-defined and is equal to the discrete 
coupling functional in~\eqref{eq:2} for all continuous $u$. 

Note that the existence and uniqueness of the solution to a boundary value 
problem with nonlocal boundary conditions is non-trivial, as one can see 
considering the following simple problem: Find~$u: [0,1] \rightarrow 
\mathbb{R}$, such that
\[
    u''= f \quad \text{in $\Omega = (0,1)$} \\
  \]
{and}
\[
    u + \mathcal{L}[u] = 0 \quad \text{on $\partial \Omega$,} \\
  \]
{with}
\[
    \mathcal{L}[u] = -\int_0^1 u \du \Omega.
\]
For any solution~$u(x)$ and constant~$\lambda$, the function~$u(x)+\lambda$ 
also solves the problem. Hence, it has no unique solution. A similar 
situation is depicted for a specific choice of the functional 
$\mathcal{L}[u]$ in \autoref{fig:10} (see Section~\ref{sec:6}). 

A full theoretical investigation of the model problem \eqref{eq:1} and its 
variational formulation is beyond the scope of this paper. In fact, the 
existence and uniqueness of the solution might be difficult to prove for such 
problems formulated on general domains. Existence and uniqueness studies are 
presented, for example, in papers~\cite{Avalishvili2011,Avalishvili2018}. In 
this paper, we assume that the model problem~\eqref{eq:1} has a unique, 
sufficiently smooth solution. 

\section{Variational formulations}

\label{sec:3}

In this section, we present several variational formulations of the model 
problem~\eqref{eq:1}. In particular, the variational formulations based on a 
strong or a weak imposition of the inhomogeneous Dirichlet-type nonlocal 
boundary condition (\hyperref[eq:1]{1b}) are given.

\subsection{Strong imposition of {Dirichlet}-type nonlocal boundary conditions}

Let us introduce the trial function space
\[
    \mathcal{S} = \{u \in H^1(\Omega): u + \mathcal{L}_1[u] = g \text{ on } \Gamma_D\}\]
and the test function space
\[
    \mathcal{V} = \{w \in H^1(\Omega): w = 0 \text{ on } \Gamma_D\}.
\]
In this way, the Dirichlet-type nonlocal boundary condition (\hyperref[eq:1]{1b}) is 
built into the space of trial functions, i.e.~this condition is imposed 
\emph{strongly}. 

After multiplying (\hyperref[eq:1]{1a}) by $w \in \mathcal{V}$ and integrating by 
parts, the identity
\begin{equation}
\label{eq:4}
    -\int_{\Omega}{\kappa \nabla w \cdot \nabla u}\du\Omega + \int_{\Gamma}{w (\kappa \nabla u \cdot \mathbf{n})}\du\Gamma + \int_{\Omega}{w f}\du\Omega = 0
\end{equation}
is obtained. Then, the Neumann-type nonlocal boundary condition 
(\hyperref[eq:1]{1c}), as well as the constraint (homogeneous boundary condition) 
$w\vert_{\Gamma_D} = 0$, is applied, and the variational formulation of problem 
\eqref{eq:1} is given as follows: Find $u \in \mathcal{S}$ such that
\begin{equation}
\label{eq:5}
    a(w, u) = \ell(w) \quad \forall w \in \mathcal{V},
\end{equation}
where
 \[
    a(w, u) = \int_{\Omega}{\kappa \nabla w \cdot \nabla u}\du\Omega + \int_{\Gamma_N}{w \mathcal{L}_2[u]}\du\Gamma \\
    \]
{and}
\[
    \ell(w) = \int_{\Omega}{w f}\du\Omega + \int_{\Gamma_N}{w h}\du\Gamma.\]

The variational formulation \eqref{eq:5} is consistent with the strong form 
\eqref{eq:1}, i.e.~the exact solution $u$ of the problem \eqref{eq:1} is also 
the solution of the variational problem \eqref{eq:5}. 

\subsection{Weak imposition of {Dirichlet}-type nonlocal boundary condition}

For a \emph{weak} imposition of the inhomogeneous Dirichlet-type nonlocal 
boundary condition (\hyperref[eq:1]{1b}), methods based on modifications of the 
identity \eqref{eq:4} can be used. 

Let us assume now that the trial and test function space both are 
$H^1(\Omega)$, i.e.~$\mathcal{S} = \mathcal{V} = H^1(\Omega)$. That is, the 
inhomogeneous Dirichlet-type nonlocal condition is no longer enforced in the 
trial function space nor the homogeneous classical Dirichlet condition is 
enforced in the test function space. Instead, the identity \eqref{eq:4} is 
modified by augmenting it with two additional terms:
\[
\begin{split}
    & -\int_{\Omega}{\kappa \nabla w \cdot \nabla u}\du\Omega + \int_{\Gamma}{w (\kappa \nabla u \cdot \mathbf{n})}\du\Gamma+ \int_{\Omega}{w f}\du\Omega \\
    & \quad + \alpha\int_{\Gamma_D}{(\kappa \nabla w \cdot \mathbf{n})(u+\mathcal{L}_1[u]-g)}\du\Gamma - \beta\int_{\Gamma_D}{w(u+\mathcal{L}_1[u]-g)}\du\Gamma = 0,
\end{split}\]
where $\beta$ is the \emph{penalty parameter} (a positive constant), and the 
constant $\alpha = 0$ or $1$. If $\alpha = 0$, we have the so-called 
\emph{penalty method}, while the method with $\alpha = 1$ is referred to as 
\emph{Nitsche's method}. The variational formulation of the problem 
\eqref{eq:1} now is stated as follows: Find $u \in H^1(\Omega)$ such that
\begin{equation}
\label{eq:6}
    a(w, u) = \ell(w) \quad \forall w \in H^1(\Omega),
\end{equation}
where
\[
\begin{split}
    a(w, u) &= \int_{\Omega}{\kappa \nabla w \cdot \nabla u}\du\Omega + \int_{\Gamma_N}{w \mathcal{L}_2[u]}\du\Gamma - \int_{\Gamma_D}{w (\kappa \nabla u \cdot \mathbf{n})}\du\Gamma \\
    & \quad - \alpha\int_{\Gamma_D}{(\kappa \nabla w \cdot \mathbf{n})(u+\mathcal{L}_1[u])}\du\Gamma + \beta\int_{\Gamma_D}{w(u+\mathcal{L}_1[u])}\du\Gamma
\end{split}
   \]
{and}\[
    \ell(w) = \int_{\Omega}{w f}\du\Omega + \int_{\Gamma_N}{w h}\du\Gamma - \alpha\int_{\Gamma_D}{(\kappa \nabla w \cdot \mathbf{n})g}\du\Gamma + \beta\int_{\Gamma_D}{w g}\du\Gamma.
\]

The variational formulation \eqref{eq:6} is consistent with the strong 
formulation of the problem \eqref{eq:1}. 

In case of classical boundary conditions ($\mathcal{L}_1 \equiv 0$, 
$\mathcal{L}_2 \equiv 0$), $\alpha = 1$ ensures the symmetry of the bilinear 
form $a(w, u)$. In case of nonlocal conditions, the bilinear form $a(w, u)$, 
in general, is non-symmetric. Asymmetry is a typical property of nonlocal 
problems and their discretisations. However, as we will see from the results 
obtained by our numerical study, in comparison to the penalty method ($\alpha 
= 0$), the terms with $\alpha = 1$ (Nitsche's method) can slightly increase 
the accuracy of the results. 

\section{Isogeometric analysis: a brief overview}

\label{sec:4}

The numerical experiments will be performed in the isogeometric framework. 
Therefore, a brief overview of the concepts and ideas used in isogeometric 
analysis is given in this section. First of all, B-splines and NURBS, which 
are the basic computer-aided geometric design (CAGD) tools applied in IGA, 
are introduced. Then, the isoparametric approach for the approximation of an 
unknown solution is shortly described. A comprehensive description of NURBS, 
their properties and related algorithms can be found, for example, 
in~\cite{Piegl1997,Rogers2001}. For more details related to IGA, 
see~\cite{Hughes2005}. 

\subsection{B-splines and {{NURBS}}}

First of all, let us define B-splines, which are progenitors of NURBS. A 
finite sequence of non-decreasing real values $\Xi = (\xi_1, \xi_2, \dotsc, 
\xi_{n+p+1})$, with $\xi_i < \xi_{i+p}$ for all $i\in\{2,\ldots,n\}$ and, 
without loss of generality, $\xi_1=0$ and $\xi_{n+p+1}=1$, is called a 
\emph{knot vector}. Here $p$ is the polynomial order (degree) of the B-spline 
basis functions and $n$ is the number of basis functions. A knot vector in 
which the first and last knot values are repeated $p + 1$ times is called an 
\emph{open} knot vector. Such knot vectors satisfy simple interpolation 
properties for function value and derivatives up to order $p$ and therefore 
play an important role both in CAGD and IGA. For a given knot vector $\Xi$, 
the B-spline basis functions of degree zero are defined as
\[
    N_{i, 0}(\xi) =
    \left\{
\begin{aligned}
            & 1 && \text{if $\xi_i \leqslant \xi < \xi_{i+1}$,} & \\
            & 0 && \text{otherwise.} &
\end{aligned}
    \right.
\]
For $p > 0$, the B-spline basis functions $N_{i, p}(\xi)$ can be constructed 
by the \emph{Cox--de Boor recursive formula}
\[
    N_{i, p}(\xi) = \frac{\xi-\xi_i}{\xi_{i+p}-\xi_i}N_{i, p-1}(\xi) + \frac{\xi_{i+p+1}-\xi}{\xi_{i+p+1}-\xi_{i+1}}N_{i+1, p-1}(\xi),\]
with the convention $\dfrac{0}{0} = 0$.

Let $\Xi = (\xi_1, \xi_2, \dotsc, \xi_{n+p+1})$ and, in the two-dimensional 
case, $\mathcal{H} = (\eta_1, \eta_2, \dotsc, \eta_{m+q+1})$ be the knot 
vectors, where $p$ and $q$ are polynomial orders of B-splines, and $n$ and 
$m$ are numbers of basis functions. Then, the univariate and bivariate NURBS 
basis functions are defined as follows:
\[
    R_i^p(\xi) = \frac{\omega_i N_{i, p}(\xi)}{\sum_{\hat{i}=1}^{n}{\omega_{\hat{i}} N_{\hat{i}, p}(\xi)}} \\
  \]
{and}
\[
    R_{i, j}^{p, q}(\xi, \eta) = \frac{\omega_{i, j} N_{i, p}(\xi) M_{j, q}(\eta)}{\sum_{\hat{i}=1}^{n}{\sum_{\hat{j}=1}^{m}{\omega_{\hat{i}, \hat{j}} N_{\hat{i}, p}(\xi) M_{\hat{j}, q}(\eta)}}},\]
where $\omega_i$ and $\omega_{i, j}$ are \emph{weights} (real positive 
numbers), $N_{i, p}$ and $M_{i, q}$ are the B-spline basis functions defined 
over the knot vectors $\Xi$ and $\mathcal{H}$, respectively. The NURBS curves 
and surfaces are constructed as linear combinations of NURBS basis functions: 
\[
    \mathbf{C}(\xi) = \sum_{i=1}^{n}{\mathbf{B}_i R_i^p(\xi)} \\
  \]
{and}
\[
    \mathbf{S}(\xi, \eta) = \sum_{i=1}^{n}{\sum_{j=1}^{m}{\mathbf{B}_{i, j} R_{i, j}^{p, q}(\xi, \eta)}},\]
where $\mathbf{B}_i$ and $\mathbf{B}_{i, j} \in \mathbb{R}^d$ ($d = 2, 3$) 
are the \emph{control points}. Since B-splines can be interpreted as NURBS 
where all weights are equal, we use the notation $N(\boldsymbol{\xi})$ to 
refer to both B-splines and NURBS.
 
In order to impose the Dirichlet-type nonlocal boundary condition 
(\hyperref[eq:1]{1b}) strongly, we use \emph{Greville abscissae}~\cite{deBoor2001}, 
which are well-known in the CAGD literature and which are widely used in 
isogeometric collocation methods~\cite{Auricchio2010,Reali2015b}. For a given 
knot vector $\Xi$, the Greville abscissae are $n$ points defined as 
\[
    \overline{\xi}_i = \frac{1}{p}(\xi_{i+1}+\xi_{i+2}+\cdots+\xi_{i+p}).
\]
The Greville abscissae naturally correspond to basis functions. The basis 
function $N_{i, p}$ attains its maximum close to the Greville abscissa 
$\overline{\xi}_i$, which is given by the average of the inner knots that 
define $N_{i, p}$. The Greville abscissae corresponding to the knot vector 
$\mathcal{H}$ are defined in the same way. Consequently, one can define 
Greville abscissae $(\overline{\xi}_i,\overline{\eta}_j)$ corresponding to 
tensor-product basis functions $N_{i,j}^{p,q}$. When using open knot vectors, 
the Greville abscissae on the boundary of the domain belong to those basis 
functions, that have non-vanishing support on the boundary. 

\subsection{Isogeometric approximation}

Let us assume that the \emph{physical domain} (model geometry) $\Omega 
\subset \mathbb{R}^2$ is represented by the geometrical mapping 
\[
    \mathbf{x}: \widehat{\Omega} \rightarrow \Omega,\]
where $\widehat{\Omega} \subset \mathbb{R}^2$ is the \emph{parametric domain} 
(in our setting we always assume $\widehat{\Omega}=(0,1)^2$). We also assume 
that this mapping is invertible, i.e.~there exists the inverse mapping 
\[
    \mathbf{x}^{-1}: \Omega \rightarrow \widehat{\Omega},\]
in short, we write $\boldsymbol{\xi}(\mathbf{x})$ for the inverse. In terms 
of the basis functions, the geometrical mapping is given as \[
    \mathbf{x}(\boldsymbol{\xi}) = \sum_{A \in \boldsymbol{\eta}^s}{\mathbf{x}_A N_A(\boldsymbol{\xi})},\]
where $\mathbf{x}_A$ are the \emph{control points} and $\boldsymbol{\eta}^s$ 
is the index set.

Note that in isogeometric analysis the physical domain is often represented 
by a collection of subdomains, so-called patches. These patches are often 
trimmed, i.e.,~their parametric domain is not a full rectangle (or cuboid) 
but only a part of it, which is usually implicitly defined. We would like to 
point out, that the imposition of boundary conditions on trimmed domains is 
usually also performed weakly, see~\cite{Marussig2018}. Similarly, coupling 
between patches may be performed weakly.

The approximation of the unknown solution $u$ in the parametric domain 
$\widehat{\Omega}$ is given by
\[
    \hat{u}^h(\boldsymbol{\xi}) = \sum_{A \in \boldsymbol{\eta}^s}{u_A N_A(\boldsymbol{\xi})},\]
where $u_A$ are the \emph{control variables}. In the physical domain 
$\Omega$, the approximation of $u$ is expressed by using the 
\emph{push-forward}:
\[
    u^h(\mathbf{x}) = \hat{u}^h \circ \mathbf{x}^{-1} = \sum_{A \in \boldsymbol{\eta}^s}{u_A N_A(\boldsymbol{\xi}(\mathbf{x}))}.
\]

The control variables $\mathbf{u} = \{u_A\}$ are determined by solving the 
linear system
\[
    \mathbf{K}\mathbf{u} = \mathbf{f},\]
where $\mathbf{K}$ is the \emph{global left-hand side matrix}, $\mathbf{f}$ 
is the \emph{global right-hand side vector}. The matrix $\mathbf{K}$ and 
vector $\mathbf{f}$ can be assembled using either a Galerkin or a collocation 
method. The details of the configurations we consider are presented in the 
following section.  

\section{Discretisations of variational formulations}

\label{sec:5}

In this section, the variational formulations presented in Section~\ref{sec:3} 
are discretised. The expressions of the global left-hand side matrices and 
the right-hand side vectors of the resulting discrete linear systems are 
given.

Let us assume that $\mathcal{S}^h$ and $\mathcal{V}^h$ are the discrete 
counterparts of the trial function space $\mathcal{S}$ and the test function 
space $\mathcal{V}$ defined in Section~\ref{sec:3} ($\mathcal{S}^h \subset 
\mathcal{S}$, $\mathcal{V}^h \subset \mathcal{V}$). The approximate solutions 
of the considered variational problems \eqref{eq:5} and \eqref{eq:6} are 
represented as
\[
    u^h = v^h + g^h,\]
where
\[
    v^h = \sum_{A \in \boldsymbol{\eta}^s \setminus \boldsymbol{\eta}_g^s}{u_A N_A(\mathbf{x})},\]
and $g^h = \Pi^h{g}$ is the interpolant or projector of the function $g$ onto 
the space of basis functions that are supported on the boundary $\Gamma_D$ 
and satisfy a nonlocal constraint (\hyperref[eq:1]{1b}). The interpolant (projector) 
is expressed as
\[
    g^h = \sum_{A \in \boldsymbol{\eta}_g^s}{g_A N_A(\mathbf{x})},\]
where $\boldsymbol{\eta}_g^s$ is the set containing all indices of the basis 
functions supported on the boundary $\Gamma_D$, corresponding to Greville 
abscissae on the boundary.
 
In the IGA framework, the interpolant $\Pi^h{g}$ can be constructed by using 
the interpolation at the \emph{mapped Greville abscissae} in the physical 
domain, the images of the Greville abscissae $\overline{\boldsymbol{\xi}}$ 
under the geometry mapping. The interpolation conditions
\[
    \Pi^h{g}(\mathbf{x}(\overline{\boldsymbol{\xi}})) = g(\mathbf{x}(\overline{\boldsymbol{\xi}}))\]
must be satisfied at the Greville abscissae $\overline{\boldsymbol{\xi}}$ 
associated with the Dirichlet boundary $\Gamma_D$.
 
Another approach to obtain $\Pi^h{g}$ is based on the computation of 
$L^2$-\emph{projection}:
\[
    \int_{\Gamma_D}{w^h\bigl(\Pi^h{g}-g\bigr)}\du\Gamma = 0 \quad \forall w^h \in \mathcal{V}^h.
\]

When the Dirichlet-type nonlocal boundary condition (\hyperref[eq:1]{1b}) is imposed 
strongly, the control variables in the expressions of $v^h$ and $g^h$ are 
determined by the linear system
\[
    \underbrace{\left[
\arraycolsep=2pt\begin{array}{lr}
      \mathbf{K}^{(1)} & \mathbf{0} \\
      \;\;\;\;{\mathbf{K}^{(2)}} 
\end{array}\right]
}_{\eqcolon \mathbf{K}_{\mathrm{s}}}
    \left\{
\begin{array}{c}
        \mathbf{u}_f \\
        \mathbf{u}_g 
\end{array}
    \right\}
    =
    \left\{
\begin{array}{c}
        \mathbf{f} \\
        \mathbf{g} 
\end{array}
    \right\},\]
where
\begin{equation}
\label{eq:7}
\begin{aligned}
    & \mathbf{u}_f = \{u_A\}, \quad A \in \boldsymbol{\eta}^s \setminus \boldsymbol{\eta}_g^s, \\
    & \mathbf{u}_g = \{g_A\}, \quad A \in \boldsymbol{\eta}_g^s, \\
    & \mathbf{K}^{(1)} = [K_{AB}^{(1)}], \quad K_{AB}^{(1)} = \int_{\Omega}{\kappa \nabla N_A \cdot \nabla N_B}\du\Omega + \int_{\Gamma_N}{N_A \mathcal{L}_2[N_B]}\du\Gamma, \quad A, B \in \boldsymbol{\eta}^s \setminus \boldsymbol{\eta}_g^s, \\
    & \mathbf{f} = \{f_A\}, \quad f_A = \int_{\Omega}{N_A f}\du\Omega + \int_{\Gamma_N}{N_A h}\du\Gamma, \quad A \in \boldsymbol{\eta}^s \setminus \boldsymbol{\eta}_g^s.
\end{aligned}
\end{equation}
However, the structure of the matrix $\mathbf{K}^{(2)}$ and the vector 
$\mathbf{g}$ depends on the method which is used for the imposition of the 
inhomogeneous Dirichlet-type nonlocal boundary condition. In case of the 
interpolation at the mapped Greville abscissae, $\mathbf{K}^{(2)}$ and 
$\mathbf{g}$ are expressed as
\[
\begin{aligned}
    & \mathbf{K}^{(2)} = [K_{AB}^{(2)}], \quad K_{AB}^{(2)} = N_B(\mathbf{x}(\overline{\boldsymbol{\xi}}_A))+\mathcal{L}_1[N_B](\mathbf{x}(\overline{\boldsymbol{\xi}}_A)), \quad A \in \boldsymbol{\eta}_g^s, \, B \in \boldsymbol{\eta}^s, \\
    & \mathbf{g} = \{g_A\}, \quad g_A = g(\mathbf{x}(\overline{\boldsymbol{\xi}}_A)), \quad A \in \boldsymbol{\eta}_g^s.
\end{aligned}
\]
If $L^2$-projection is applied, then
\[
\begin{aligned}
    & \mathbf{K}^{(2)} = [K_{AB}^{(2)}], \quad K_{AB}^{(2)} = \int_{\Gamma_D}{N_A(N_B+\mathcal{L}_1[N_B])}\du\Gamma, \quad A \in \boldsymbol{\eta}_g^s, \, B \in \boldsymbol{\eta}^s, \\
    & \mathbf{g} = \{g_A\}, \quad g_A = \int_{\Gamma_D}{N_A g}\du\Gamma, \quad A \in \boldsymbol{\eta}_g^s.
\end{aligned}
\]

The weak imposition of the inhomogeneous Dirichlet-type nonlocal boundary 
condition (\hyperref[eq:1]{1b}) leads to the linear system
\[
    \underbrace{(\mathbf{K}^{(1)} - \mathbf{K}^{(3)} - \alpha \mathbf{K}^{(4)} + \beta \mathbf{K}^{(5)})}_{\eqcolon \mathbf{K}_{\mathrm{w}}}\mathbf{u} = \mathbf{f} - \alpha \mathbf{g}^{(1)} + \beta \mathbf{g}^{(2)},\]
where $\mathbf{K}^{(1)}$ and $\mathbf{f}$ are defined as in \eqref{eq:7} with 
indices $A, B \in \boldsymbol{\eta}^s$, and
\begin{gather*}
    \mathbf{K}^{(3)} = [K_{AB}^{(3)}], \quad K_{AB}^{(3)} = \int_{\Gamma_D}{N_A (\kappa \nabla N_B \cdot \mathbf{n})}\du\Gamma, \\
    \mathbf{K}^{(4)} = [K_{AB}^{(4)}], \quad K_{AB}^{(4)} = \int_{\Gamma_D}{(\kappa \nabla N_A \cdot \mathbf{n})(N_B+\mathcal{L}_1[N_B])}\du\Gamma, \\
    \mathbf{K}^{(5)} = [K_{AB}^{(5)}], \quad K_{AB}^{(5)} = \int_{\Gamma_D}{N_A(N_B+\mathcal{L}_1[N_B])}\du\Gamma, \\
    \mathbf{g}^{(1)} = \{g_A^{(1)}\}, \quad g_A^{(1)} = \int_{\Gamma_D}{(\kappa \nabla N_A \cdot \mathbf{n})g}\du\Gamma, \\
    \mathbf{g}^{(2)} = \{g_A^{(2)}\}, \quad g_A^{(2)} = \int_{\Gamma_D}{N_A g}\du\Gamma, \\
    A, B \in \boldsymbol{\eta}^s.
\end{gather*}

In contrast to the standard practice in the finite element analysis using 
models involving classical boundary conditions, the boundary control 
variables $\mathbf{u}_g$ are solved at the same time as the interior control 
variables $\mathbf{u}_f$.
 
Due to the nonlocal boundary conditions, the global left-hand side matrices 
$\mathbf{K}_{\mathrm{s}}$ and $\mathbf{K}_{\mathrm{w}}$, typically, are 
non-symmetric. Some examples of their sparsity patterns will be presented in 
the next section.

\section{Numerical study}

\label{sec:6}

In this section we present the results of the numerical study which was 
conducted in order to examine and compare the methods for the imposition of 
inhomogeneous Dirichlet-type nonlocal boundary conditions.
 
\subsection{Implementation details}

The approximate solution methods for the model problem \eqref{eq:1} have been 
implemented using \textsc{GeoPDEs} (version 3.1.0), a free IGA software 
suite~\cite{Vazquez2016,GeoPDEs2016}. \textsc{GeoPDEs} can be launched both 
in \textsc{MATLAB} and \textsc{GNU Octave} environments.
 
\begin{figure}
\centering
\includegraphics[scale=1.0]{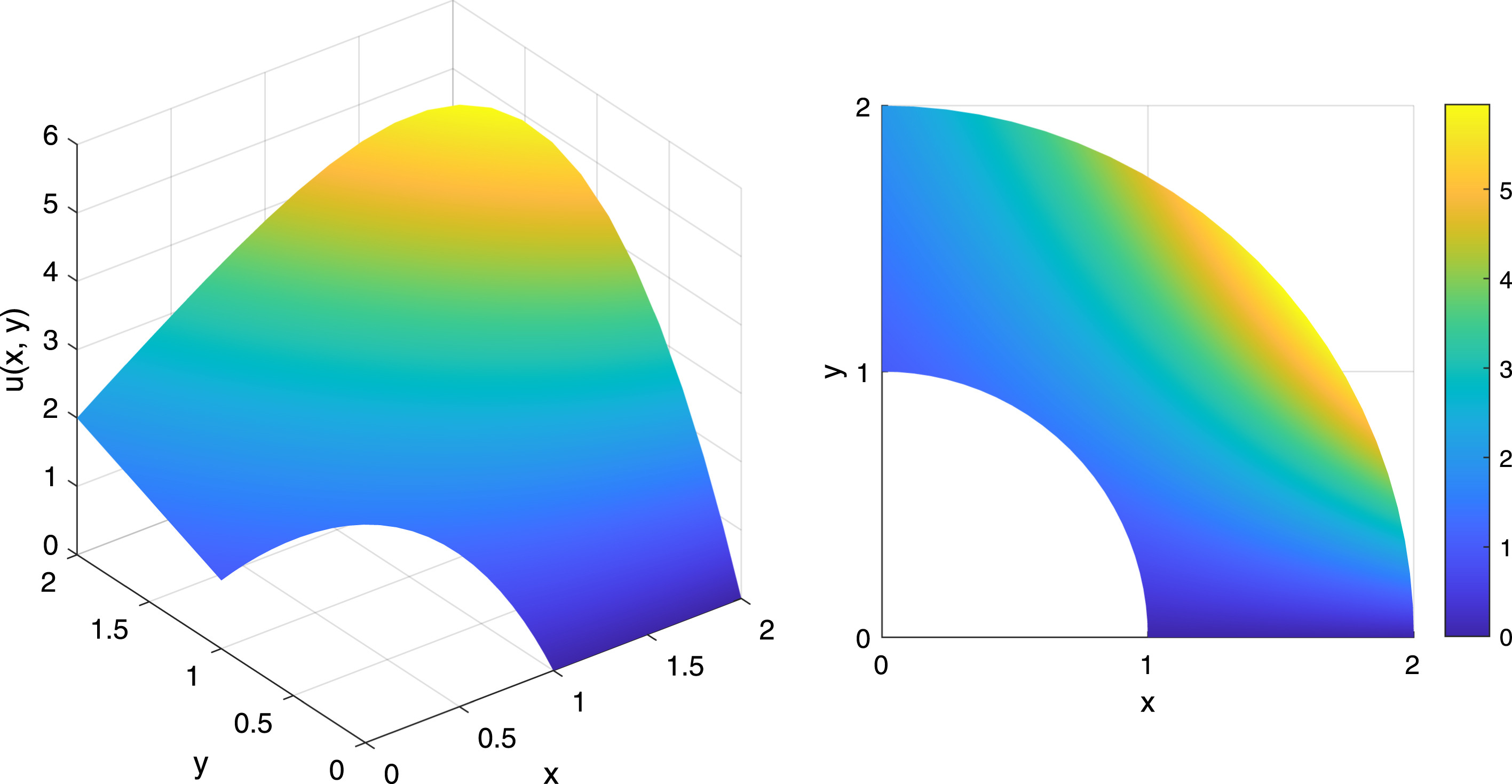}
\caption{The manufactured solution of the test problems. \coltext
 \label{fig:1}}
\end{figure}

\begin{figure}
\centering
\includegraphics[scale=1.00]{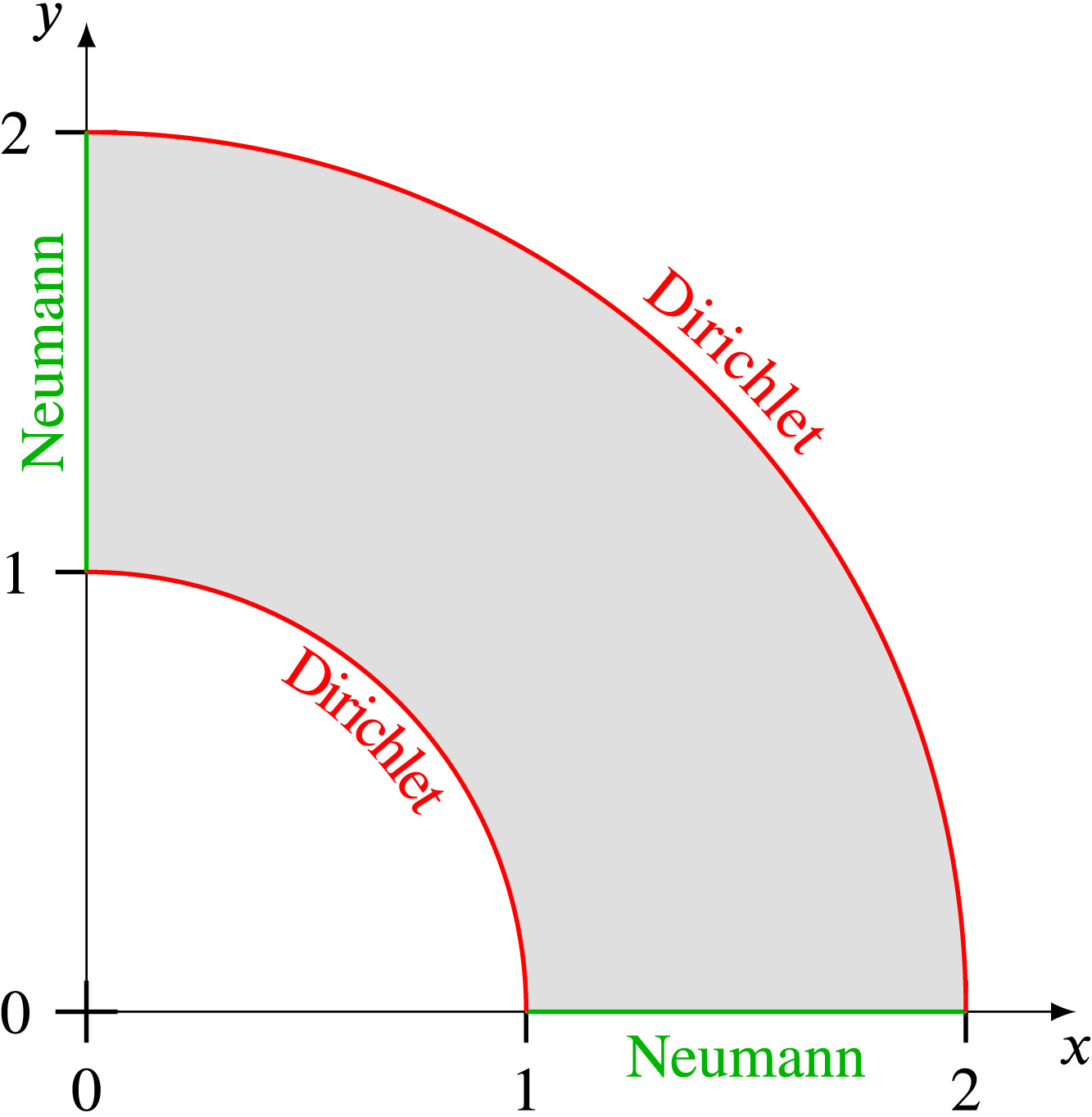}
\caption{The description of the test problems.\label{fig:2}}
\end{figure}

The two-dimensional ($d = 2$) problems with different (discrete or integral) 
nonlocal boundary conditions and the manufactured solution
\[
    u(x, y) = \exp{(x)} \cdot y,
\]
 see \autoref{fig:1}, were analysed. The diffusivity coefficient in 
(\hyperref[eq:1]{1a}) was assumed to be $\kappa(x, y) \equiv 1$.   The problems were 
formulated on a quarter of a ring with the inner radius being equal to $1$ 
and the outer radius equal to $2$ (see \autoref{fig:2}). The boundary part 
for the Neumann-type nonlocal boundary condition (\hyperref[eq:1]{1b}) was assumed to 
be
\[
    \Gamma_N = \{(x, y): 1 \leqslant x \leqslant 2, y = 0\} \cup \{(x, y): x = 0, 1 \leqslant y \leqslant 2\},\]
while the Dirichlet-type nonlocal boundary condition (\hyperref[eq:1]{1b}) was 
imposed on the rest of the boundary $\Gamma$ ($\Gamma_D = \Gamma \setminus 
\Gamma_N$). The domain parameterisation distributed together with 
\textsc{GeoPDEs} toolbox~\cite{GeoPDEs2016} in the file \texttt{ring.mat} was 
used.
 
\begin{figure}
\centering
\includegraphics[scale=1.0]{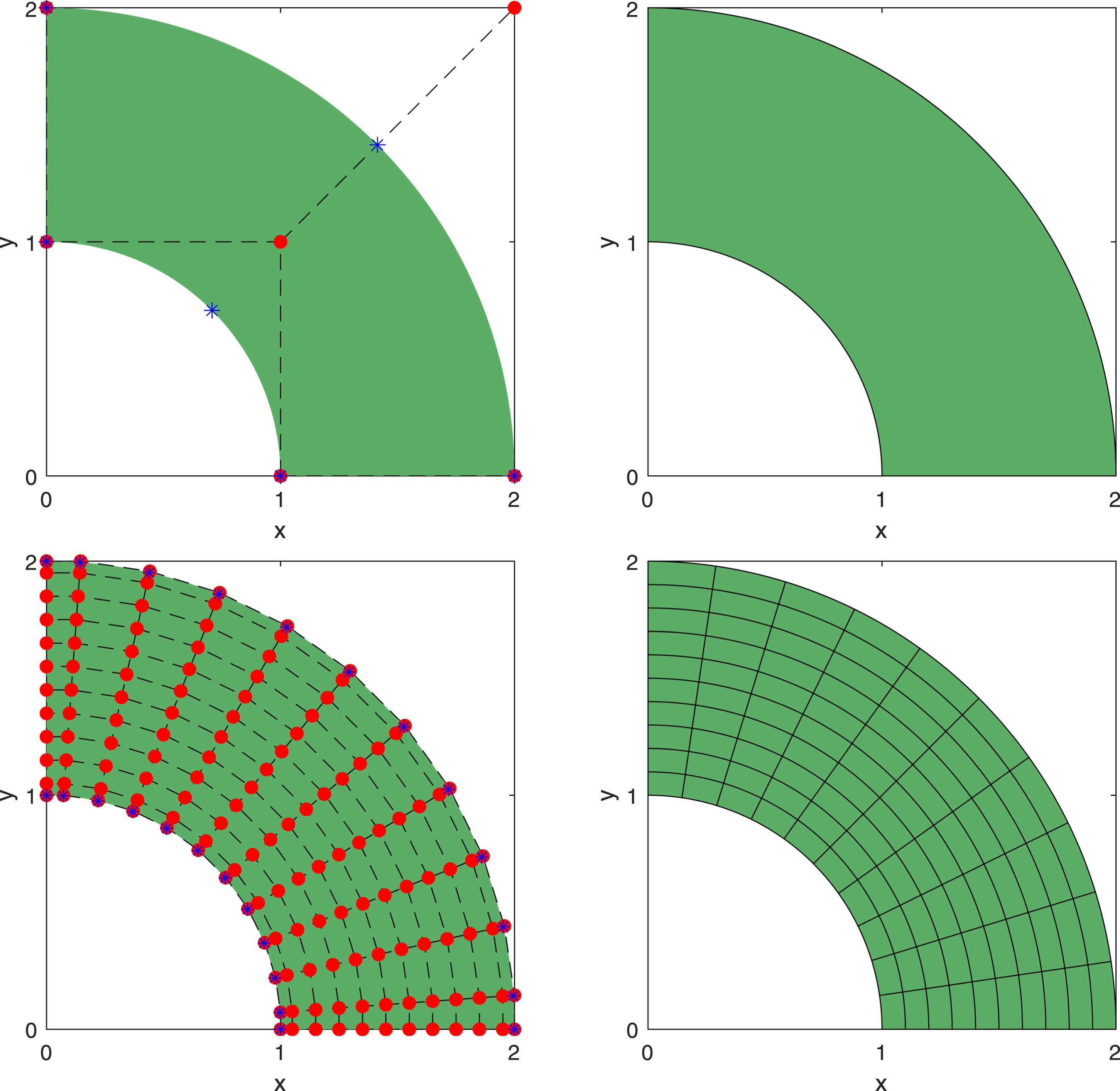}
\caption{The control meshes (left) and the physical meshes (right) before 
  (top) and after (bottom) $k$-refinement. The control points are denoted 
  as dots (red) and the images of the Greville abscissae associated with 
  the boundary $\Gamma_D$ are denoted as stars (blue). \coltext
\label{fig:3}}
\end{figure}

Two cases of nonlocal boundary conditions (\hyperref[eq:1]{1b}) and (\hyperref[eq:1]{1c}) 
were considered:
\begin{description}
  \item[Case~1 (Discrete conditions).] $\mathcal{L}_1[u] = \mathcal{L}_2[u] = \gamma u(1, 1)$
  \item[Case~2 (Integral conditions).] $ \mathcal{L}_1[u] = \mathcal{L}_2[u] = {\displaystyle \gamma \int_{\Omega}{u}\du\Omega}$
\end{description}

The role of the parameter (a real number) $\gamma$ is to control the 
influence of nonlocal terms. The case of $\gamma = 0$ corresponds to the 
classical Dirichlet and Neumann boundary conditions (classical case). It is 
assumed that $\gamma = 1$, if it is not mentioned otherwise.
 
In order to enforce a nonlocal discrete coupling which involves the value of 
the unknown function at the point $\mathbf{x}^*$ in the physical domain 
$\Omega$ (as defined by the functionals of Case~1), one needs to identify the 
corresponding point $\boldsymbol{\xi}^*$ in the parametric domain 
$\widehat{\Omega}$ by calculating the \emph{pull-back} $\boldsymbol{\xi}^* = 
\mathbf{x}^{-1}(\mathbf{x}^*)$, i.e.~by solving the equation 
$\mathbf{x}({\boldsymbol{\xi}}^*) = \mathbf{x}^*$. For this, the standard 
\textsc{MATLAB}/\textsc{GNU Octave} routine \texttt{fsolve} was used. 
Meanwhile, nonlocal integral couplings (as in Case~2) can be implemented 
using the integration techniques which are typically applied in the IGA 
framework.
 
The numerical results presented in this section were obtained using NURBS of 
polynomial orders $p = q = 2$, $C^1$ continuity in both parametric directions 
and the physical mesh consisting of {10~\ttimes~10} knot spans. The control 
and physical meshes before and after 
$k$-refinement~\cite{Hughes2005,Cottrell2009} are shown in \autoref{fig:3}. 
In all numerical experiments, the integration was performed using a standard, 
element-wise, $(p+1) \times (q+1)$ tensor-product Gauss quadrature rule. The 
penalty parameter $\beta = 10^2$ was used unless it is stated otherwise. The 
resulting linear systems, as a rule, are dense therefore iterative solvers 
might not be efficient and direct solvers might be required. In the numerical 
study, the linear systems were solved using the left division operator 
'\texttt{\textbackslash}', which is standard both in \textsc{MATLAB} and 
\textsc{GNU Octave} environments.
 
The accuracy of the results was estimated using the $L^2$-norm of the 
absolute error. The 2-norm condition number $\cond{(\bullet)} = \|  \bullet 
\| _2 \cdot \| \bullet^{-1}\| _2$ was used to estimate the conditioning of the 
global left-hand side matrices.
 
\subsection{Results and discussion}

In the numerical study we focused on the investigation of the previously 
described imposition methods for inhomogeneous Dirichlet-type nonlocal 
boundary conditions. In particular, we examine the influence of the penalty 
parameter $\beta$ as well as the influence of the nonlocal boundary 
conditions on the properties of the methods. The different methods were 
compared in terms of their accuracy and conditioning.
 
\subsubsection*{Influence of the penalty parameter}
\addcontentsline{toc}{subsubsection}{Influence of the penalty parameter}

The weak methods for the imposition of inhomogeneous Dirichlet-type nonlocal 
boundary conditions (the penalty method and Nitsche's method) depend on the 
penalty parameter $\beta$ which controls the coercivity of the corresponding 
bilinear form $a(w, u)$. Too large values of the penalty parameter $\beta$ 
yield a poorly conditioned global left-hand side matrix 
$\mathbf{K}_{\mathrm{w}}$. Therefore, it is necessary to find a suitable value 
of the penalty parameter $\beta$ which guarantees the coercivity of the 
bilinear form and leads to the desired accuracy of the numerical results.
 
We performed a parametric study to investigate the influence of the penalty 
parameter $\beta$ on the accuracy and conditioning of the methods applied to 
solve the problems with the classical and nonlocal boundary conditions.
 
\begin{figure}
\centering
\includegraphics[scale=1.0]{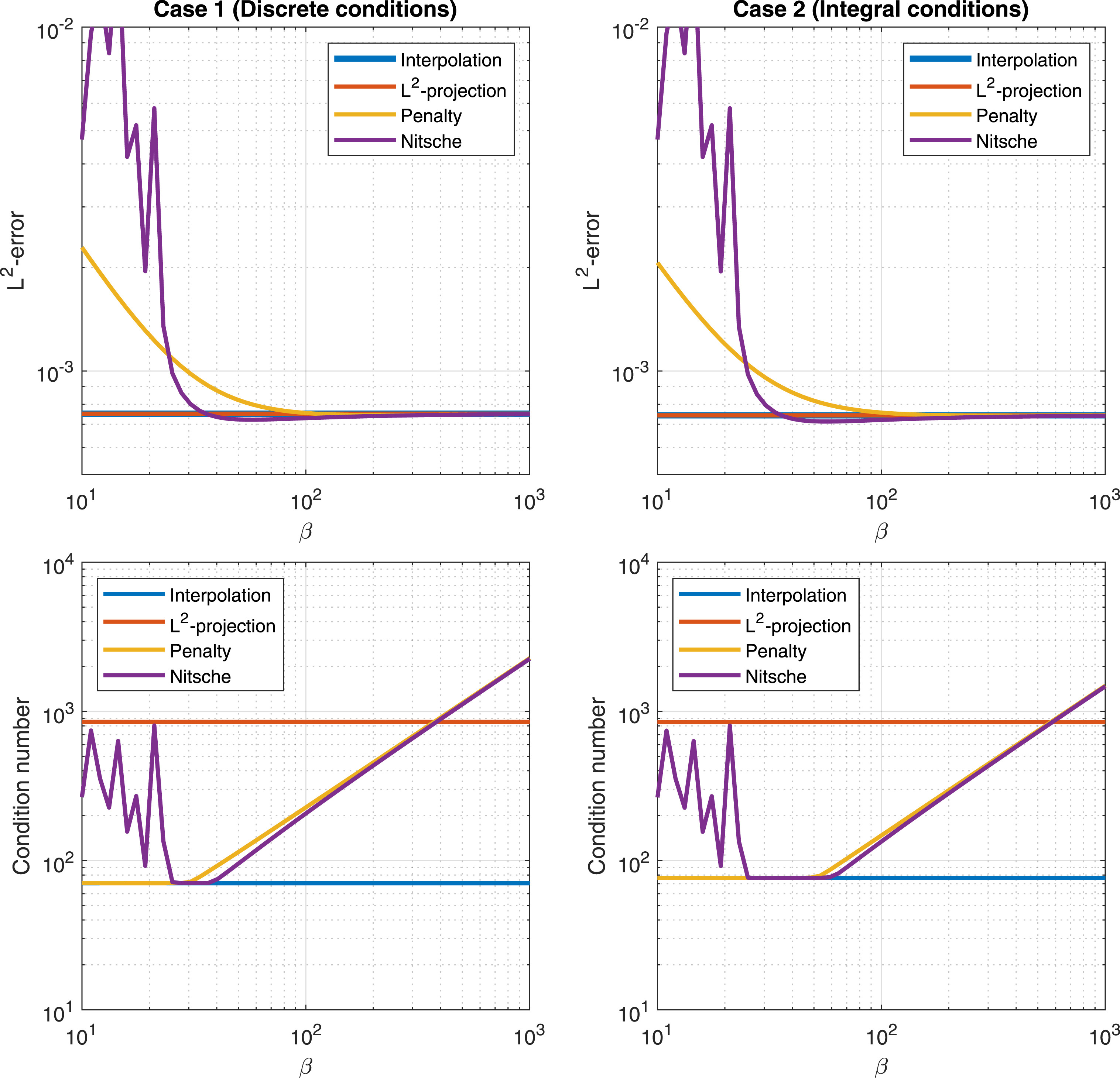}
\caption{The dependence of $L^2$-norm of the absolute error values (top) 
  and the condition numbers of the global left-hand side matrices (bottom)   
  on the values of the penalty parameter $\beta$. \coltext
\label{fig:4}}
\vskip -15pt
\end{figure}

\autoref{fig:4} depicts the $L^2$-norms of the absolute errors and the 
condition numbers of the global left-hand side matrices obtained using the 
penalty and Nitsche's methods with various values of the penalty parameter 
$\beta$. The $L^2$-norms of the absolute error values and the condition 
numbers obtained using the methods independent on the penalty parameter are 
also referenced.
 
As one can see from \autoref{fig:4} (top), the difference in the accuracy of 
the results obtained using two different methods with strongly imposed 
inhomogeneous Dirichlet-type nonlocal boundary condition (the interpolation 
at the Greville abscissae and the computation of $L^2$-projection) is 
insignificant. However, the application of the interpolation at the Greville 
abscissae leads to the methods with the global left-hand side matrices which 
has better conditioning (\autoref{fig:4}, bottom). In fact, the method based 
on the interpolation at the Greville abscissae has the best conditioning 
among all four investigated approaches.
 
Comparing the weak imposition methods, the results obtained using Nitsche's 
method are slightly more accurate than those obtained by using the penalty 
method, if the penalty parameter $\beta \sim 10^2$ (\autoref{fig:4}, top). 
For smaller values of $\beta$ however, Nitsche's method becomes more 
ill-conditioned and yields worse results. The difference becomes 
insignificant with larger values of the penalty parameter. \autoref{fig:4} 
(bottom) shows how the condition numbers of the resulting global left-hand 
side matrices grow when the value of the penalty parameter $\beta$ increases. 
It is also important to note, that Nitsche's method (which is variationally 
consistent for the problem with classical boundary conditions) does not yield 
substantially better stability properties than the penalty method when 
nonlocal boundary conditions are enforced.

\begin{figure}
\centering
\includegraphics[scale=1.0]{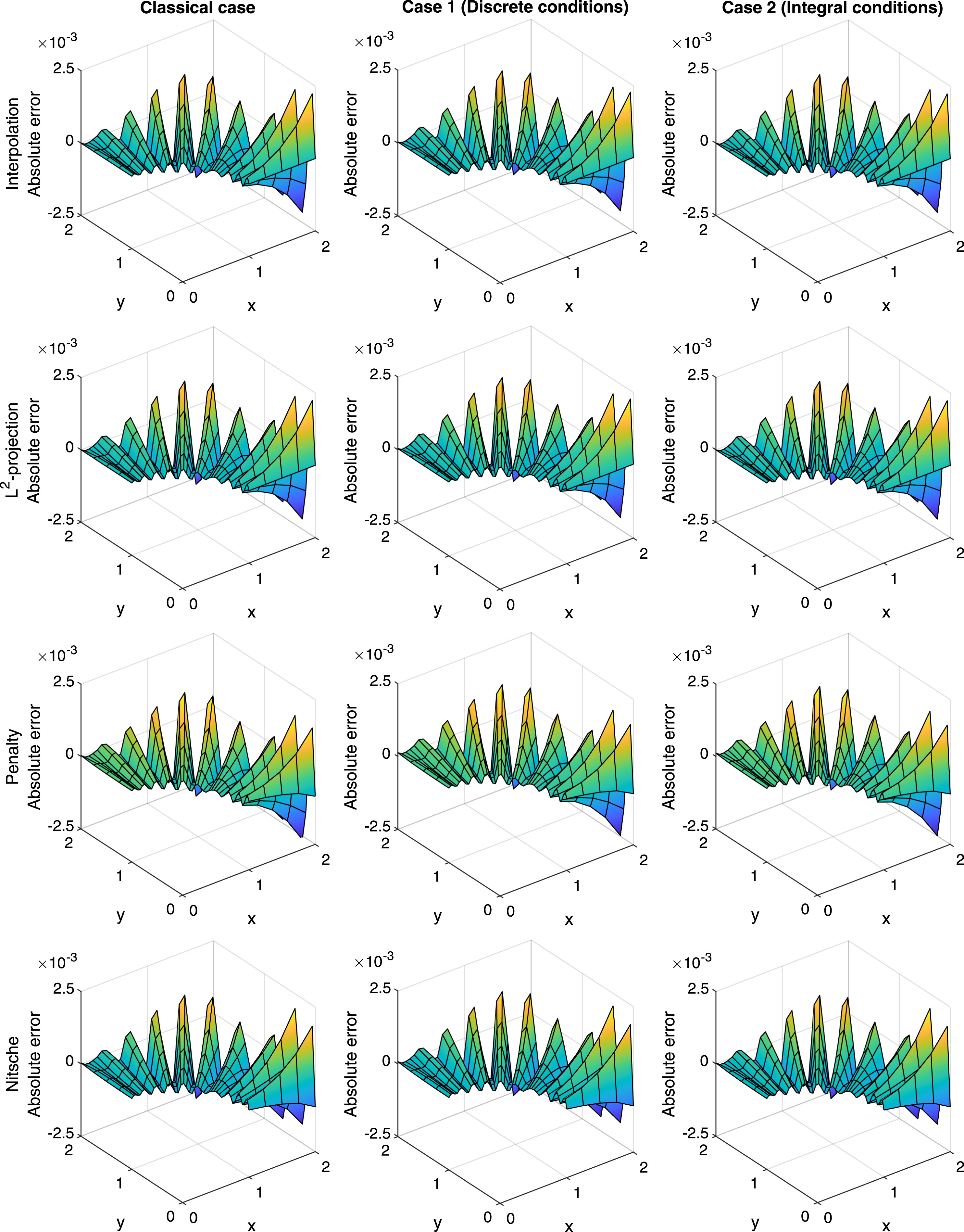}
\caption{The absolute errors of the numerical solutions obtained using 
  different methods for the implementation of the inhomogeneous 
  Dirichlet-type nonlocal boundary conditions.\label{fig:5}}
\end{figure}

In \autoref{fig:5}, we present the absolute errors of the numerical solutions 
obtained using various methods of the imposition of the inhomogeneous 
Dirichlet-type nonlocal boundary condition. Since the value of the penalty 
parameter $\beta$ used in the computations was quite large ($\beta = 10^2$), 
the numerical results obtained using the weak imposition methods in terms of 
accuracy are comparable to those obtained when the discrete or integral 
Dirichlet-type nonlocal boundary conditions are imposed strongly. Note that 
in comparison to the penalty and Nitsche's methods, the strong imposition 
methods result in lower errors at the boundary part for the Neumann-type 
boundary condition.
 
\begin{figure}
\centering
\includegraphics[scale=1.0]{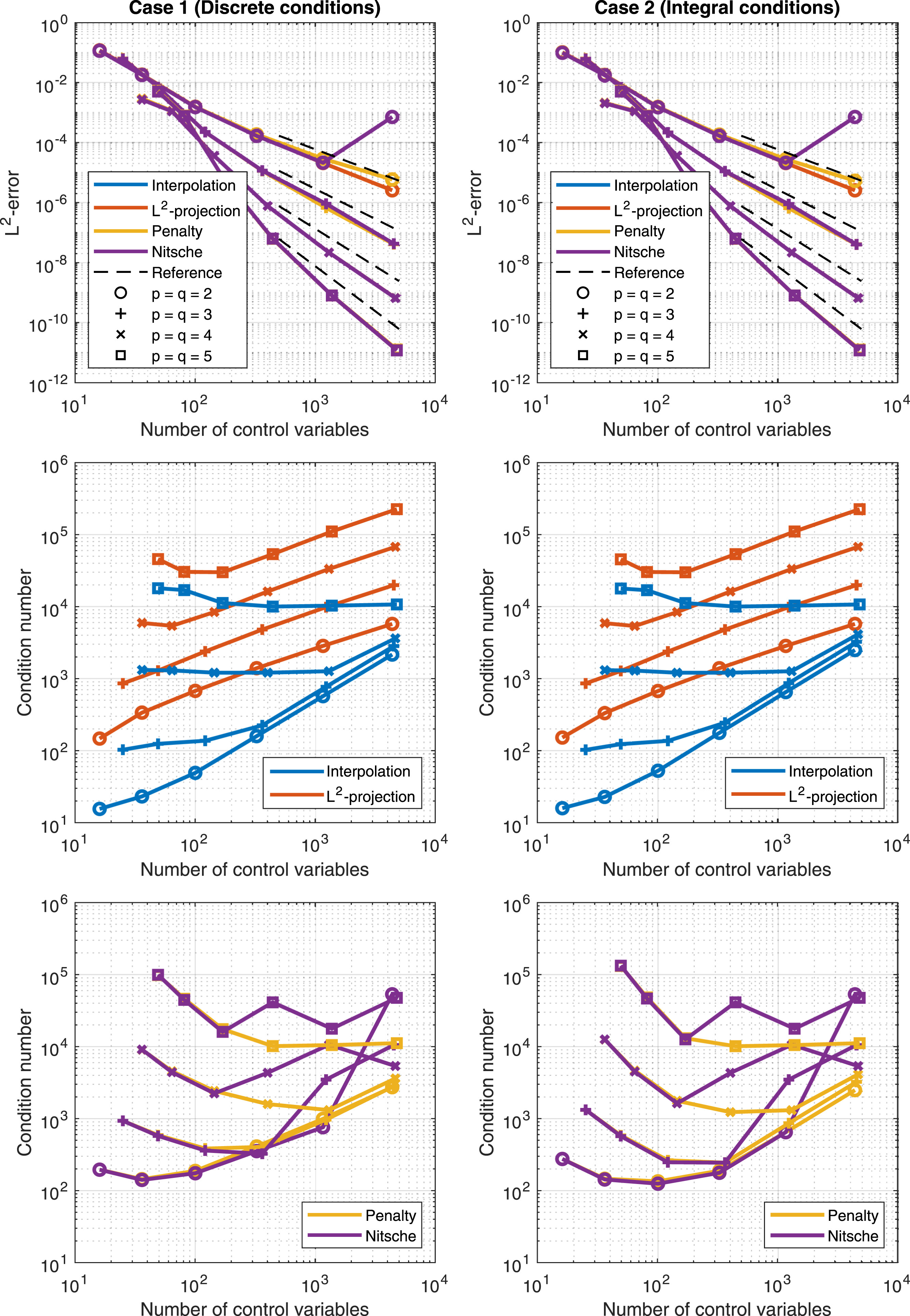}
\caption{The dependence of the $L^2$-error norms (top) and condition 
  numbers of the global left-hand side matrices (middle and bottom) on the 
  total number of control variables (degrees of freedom) under the 
  refinement using basis functions of different degree. For the weak 
  imposition methods we used a fixed penalty parameter $\beta$. Note that 
  some lines are not visible. All methods except the penalty and Nitsche's 
  method in the case $p = q = 2$ converge with the same rate. \coltext
\label{fig:06}}
\end{figure}

We investigated the convergence of the methods under knot refinement and 
degree elevation. The dependence of the $L^2$-error and of the condition 
number of the global left-hand side matrix on the number of control variables 
is presented in \autoref{fig:06}. We see that optimal convergence is achieved 
in all cases, except the case of the weak imposition methods with $p = q = 2$. 

\begin{figure}
\centering
\includegraphics[scale=1.0]{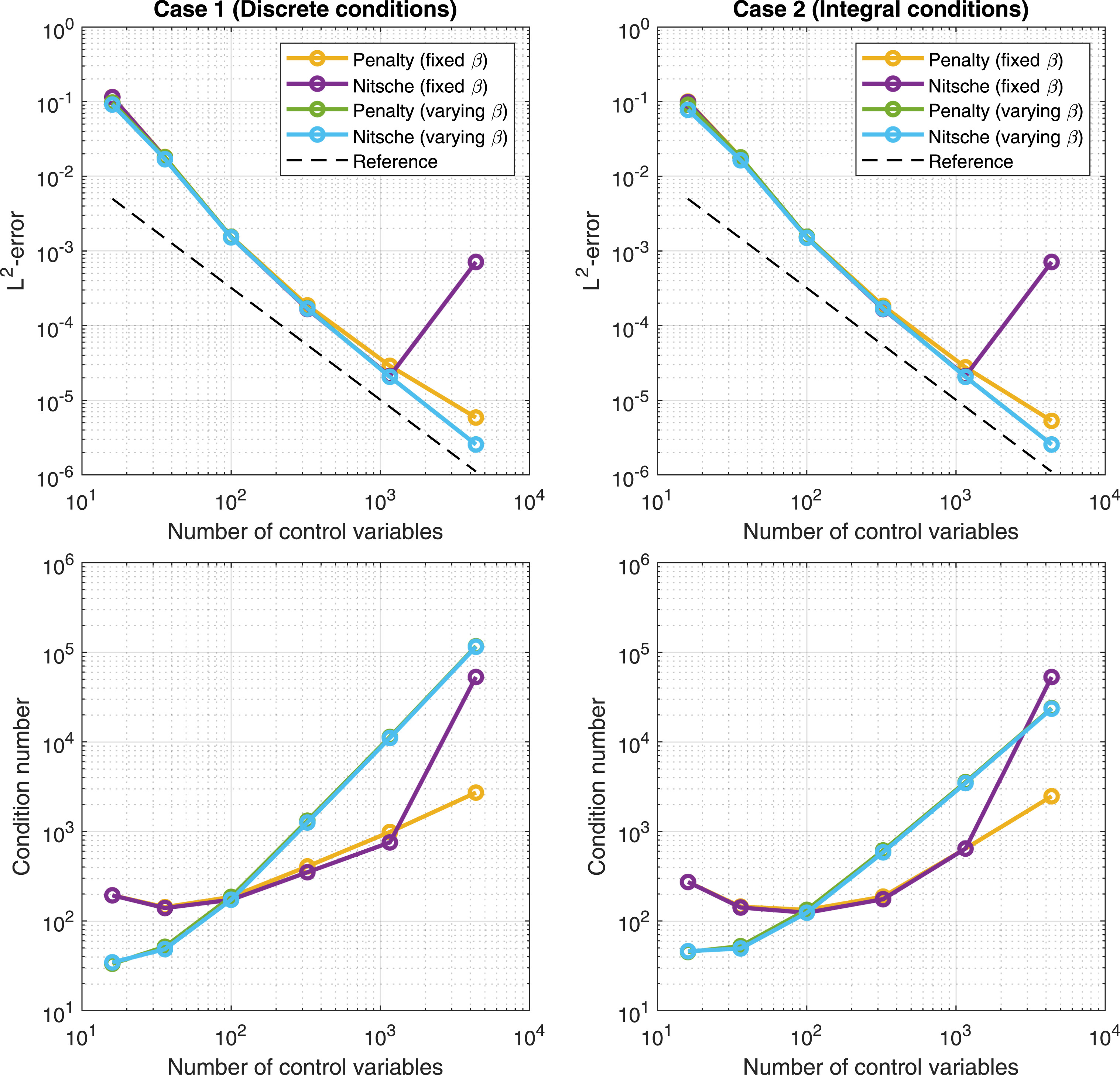}
\caption{The dependence of the $L^2$-error norms (top) and condition 
  numbers of the global left-hand side matrices (bottom) on the total number 
  of control variables (degrees of freedom). We compared a fixed and 
  varying penalty parameter $\beta$. \coltext
\label{fig:07}}
\end{figure}

Note that the results in \autoref{fig:06} have been obtained using the fixed 
penalty parameter ($\beta = 10^2$). However, the optimal convergence rate 
might be recovered if the penalty parameter value is properly adapted during 
refinement. As an example, penalty parameter values equalling the number of 
control variables ($\beta \sim 1/h^2$), following the suggestion 
from~\cite{Fernandez-Mendez2004}, were used in the experiment presented 
in \autoref{fig:07}.

\subsubsection*{Influence of nonlocal boundary conditions}
\addcontentsline{toc}{subsubsection}{Influence of nonlocal boundary conditions}

\begin{figure}
\centering
\includegraphics[scale=1.0]{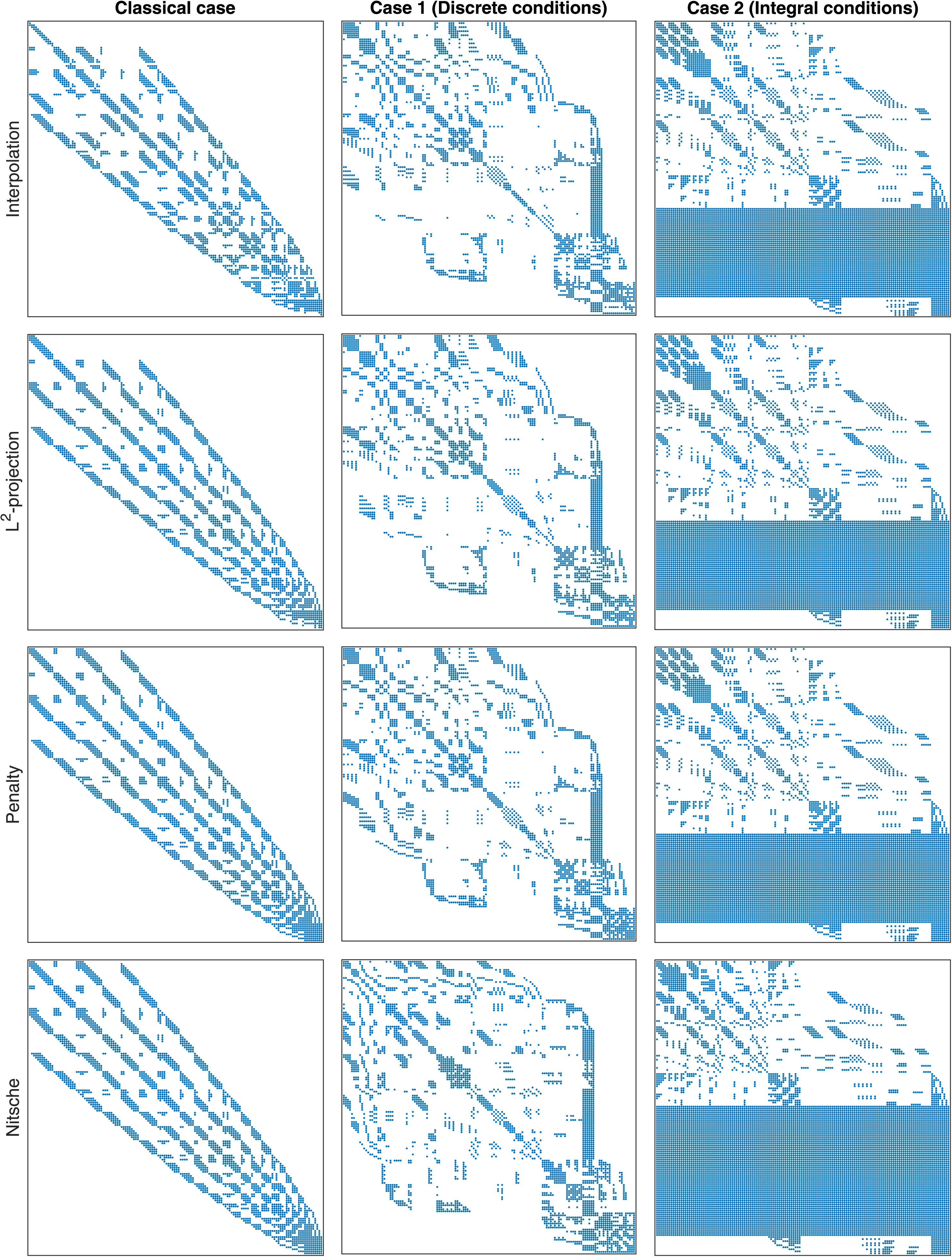}
\caption{The sparsity patterns of the symmetric reverse Cuthill--McKee
  orderings of the global left-hand side matrices.\label{fig:09}}
\end{figure}

In comparison to the classical case, the appearance of nonlocal boundary 
conditions results in additional non-zero elements in the global left-hand 
side matrices. Examples of the sparsity patterns of these matrices are 
presented in \autoref{fig:08}. Note that the structure of the global 
left-hand side matrices is implementation-dependent. In \autoref{fig:09}, the 
sparsity patterns of the symmetric reverse Cuthill--McKee orderings are 
depicted.
 
\autoref{tab:1} presents the percentage of non-zero elements in the global 
left-hand side matrices and the bandwidths of their standard and symmetric 
reverse Cuthill--McKee orderings. One can see that, in comparison to the 
classical case, nonlocal boundary conditions increase the number of non-zero 
elements. The interpolation at the Greville abscissae leads to the most 
sparse matrices, while the matrices obtained using the penalty and Nitsche's 
methods are the densest. Nonlocal integral boundary conditions increase the 
number of non-zero elements extremely and lead to full bandwidth matrices. 

A parametric study was performed with the aim to investigate how nonlocal 
conditions affect the conditioning of the global left-hand side matrices 
$\mathbf{K}_{\mathrm{s}}$ and $\mathbf{K}_{\mathrm{w}}$. The test problems were 
solved using various values of the parameter $\gamma$.

\autoref{fig:10} clearly shows that nonlocal boundary conditions have a 
negative effect on the conditioning of the global left-hand side matrices. On 
the one hand, the condition numbers increase as the absolute value of the 
parameter $\gamma$ is increased. On the other hand, the problem becomes 
ill-conditioned around special values of $\gamma$, which are visible in the 
spikes in \autoref{fig:10}. These values correspond to ill-posed problems or 
to problems which do not possess a unique solution as described in 
Section~\ref{sec:2}. The global left-hand side matrices corresponding to the 
methods with strongly imposed inhomogeneous Dirichlet-type nonlocal boundary 
conditions exhibit, in general, better conditioning in comparison to the 
matrices corresponding to the methods in which nonlocal boundary conditions 
are imposed weakly. 

\clearpage

\section{Summary and concluding remarks}

\label{sec:7}

Imposing inhomogeneous Dirichlet (essential) boundary conditions is an 
important issue in computational approaches using non-interpolatory basis 
functions. Such functions are used in various meshfree methods, as well as in 
isogeometric analysis. 

In this paper, we considered a model problem consisting of the Poisson 
equation equipped with nonlocal boundary conditions. Strong and weak methods 
for the imposition of Dirichlet-type nonlocal boundary condition were derived 
and compared experimentally in the isogeometric framework. The influence of 
nonlocal boundary conditions on the properties of the methods was 
investigated. Based on the results of numerical study, we derived the 
following conclusions:

\begin{figure}
\centering
\includegraphics[scale=1.0]{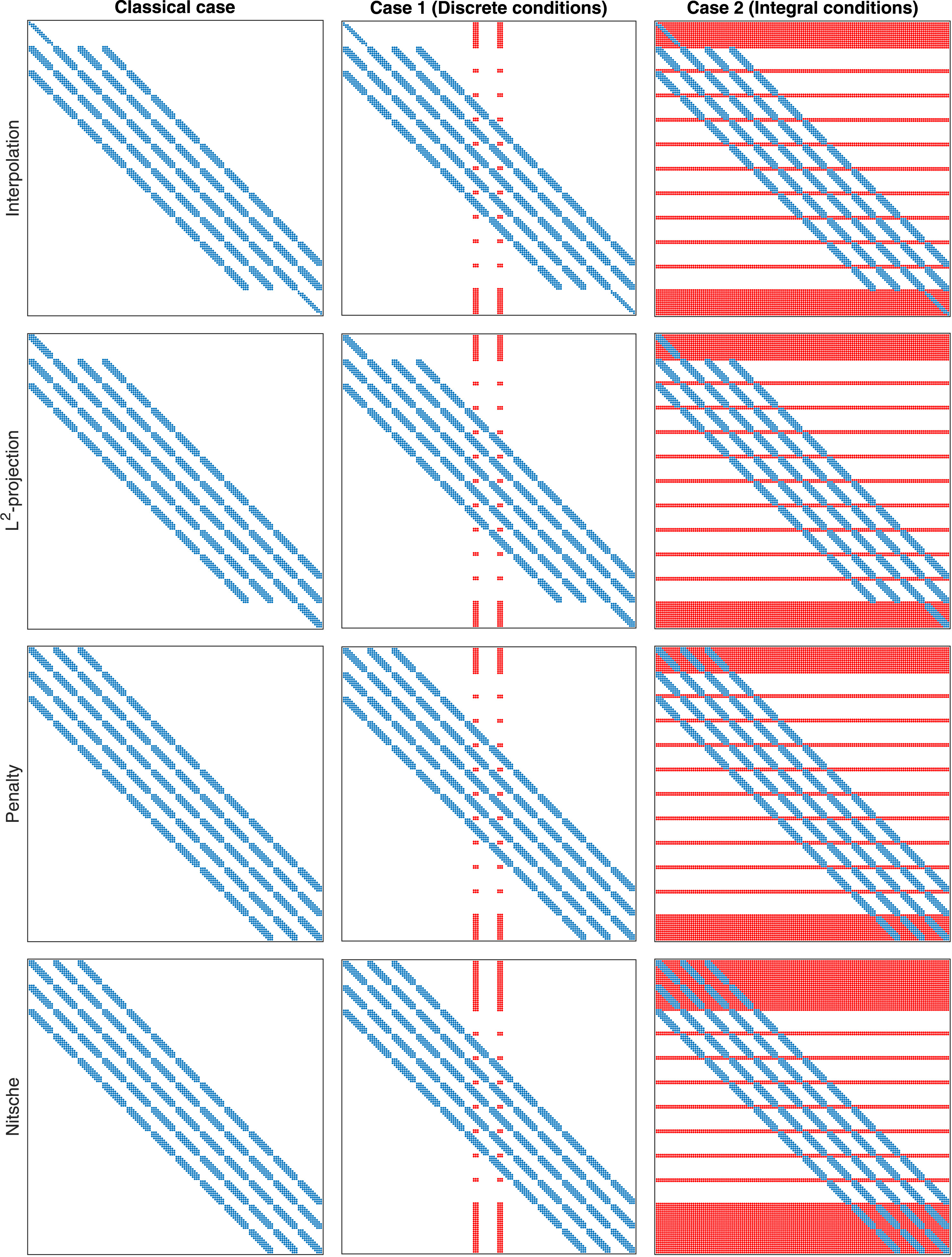}
\caption{The sparsity patterns of the global left-hand side matrices. The 
  non-zero entries arisen due to nonlocal boundary conditions are marked in 
  red. \coltext
\label{fig:08}}
\end{figure}

\clearpage
\afterpage{
\begin{figure}
\centering
\includegraphics[scale=0.99]{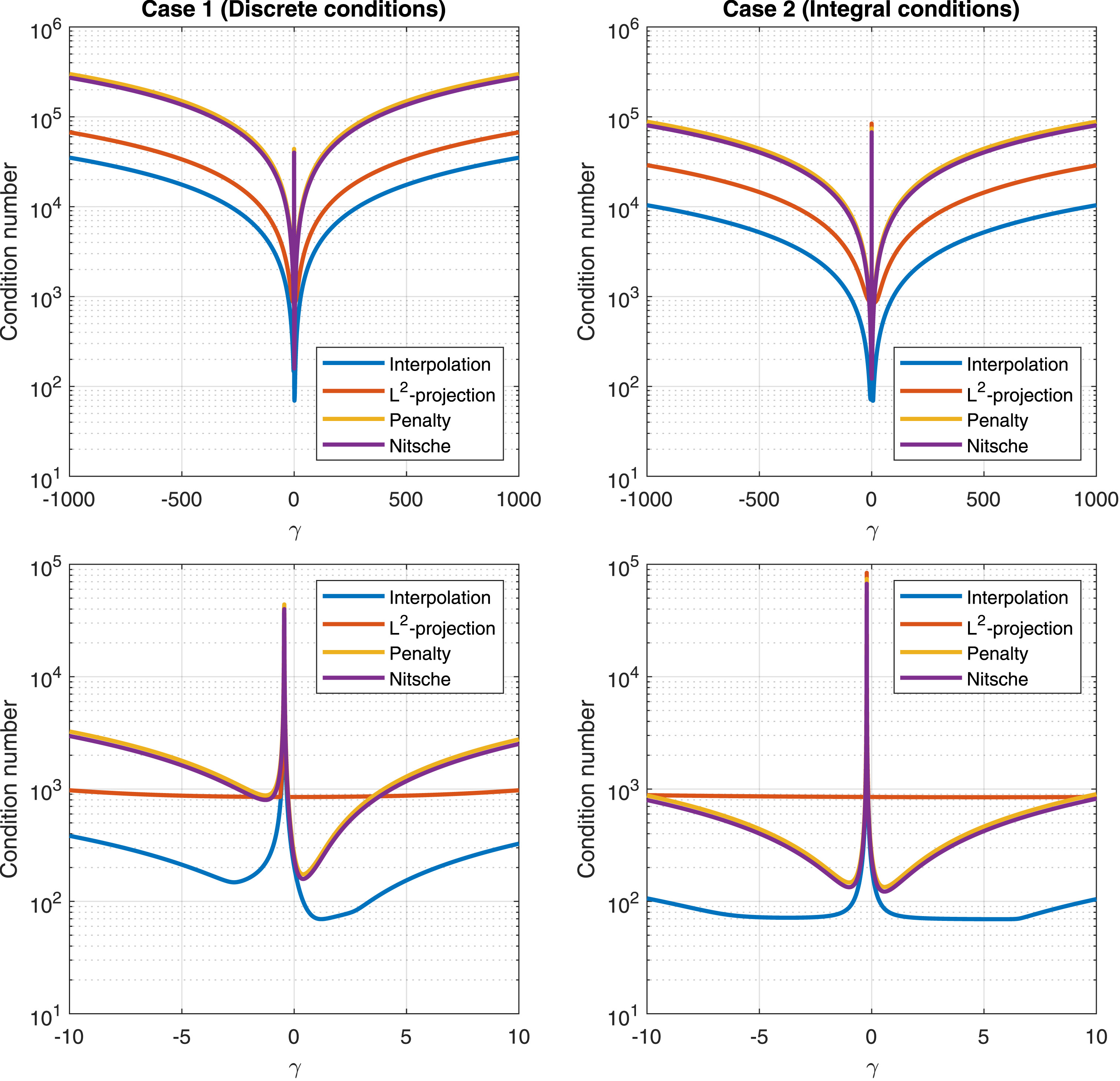}
\caption{The dependence of the condition numbers of the global left-hand 
  side matrices on the values of the parameter $\gamma$. \coltext
\label{fig:10}}
\vskip -10pt
\end{figure}

\begin{table}[width=0.75\textwidth] 
 \caption{\label{tab:1}The percentage of non-zero elements (nnz) and the 
bandwidth of the symmetric reverse Cuthill--McKee orderings of the global 
left-hand side matrices. All the matrices are of size {144~\ttimes~144}.} 
\vskip -5pt
\scriptsize
\begin{tabular}{p{3.7cm} p{3.7cm} p{3.7cm}}
  \br
                   & nnz (\%) & Bandwidth \\
  \mr
  \multicolumn{3}{l}{\emph{Classical case}} \\
  Interpolation    & 12.81    & 45 \\
  $L^2$-projection & 13.02    & 45 \\
  Penalty          & 14.06    & 45 \\
  Nitsche          & 14.06    & 45 \\[2pt]
  \multicolumn{3}{l}{\emph{Case~\emph{1} (Discrete conditions)}} \\
  Interpolation    & 14.08    & 80 \\
  $L^2$-projection & 14.29    & 80 \\
  Penalty          & 15.34    & 80 \\
  Nitsche          & 15.91    & 95 \\[2pt]
  \multicolumn{3}{l}{\emph{Case~\emph{2} (Integral conditions)}} \\
  Interpolation    & 41.67    & 135 \\
  $L^2$-projection & 41.67    & 135 \\
  Penalty          & 41.67    & 135 \\
  Nitsche          & 53.70    & 135 \\
  \br
\end{tabular}
\end{table} 
}

\clearpage
\FloatBarrier

\begin{enumerate}[\textbullet]
\itemsep=-2pt
\item Variational formulations with weakly imposed inhomogeneous 
    Dirichlet-type nonlocal boundary conditions depend strongly on the choice 
    of the penalty parameter. For small values of the penalty parameter, the 
    accuracy is reduced. For larger values of the penalty parameter, the 
    accuracy is improved, at the cost of the conditioning of the problem. The 
    strong imposition method using the interpolation at the Greville 
    abscissae has the best conditioning whereas the smallest error is reached 
    using Nitsche's method, if the penalty parameter is selected properly.
    \item The appearance of nonlocal boundary conditions leads to an     
    increased number of non-zero elements in the global left-hand side     
    matrix. The weak imposition of Dirichlet-type nonlocal boundary     
    condition leads to slightly denser global left-hand side matrices 
    than those obtained when using the strong imposition methods. The     
    global left-hand side matrices obtained when solving the problem with     
    the nonlocal integral condition are much denser in comparison to the     
    problems with nonlocal discrete boundary conditions.
    \item Numerical evidence indicates that nonlocal boundary conditions 
    might have a negative influence on the conditioning of the global 
    left-hand side matrices.
\end{enumerate}

Note that the density and structure of the resulting global left-hand side 
matrices also has an influence on the preferred choice of the solver, which 
may be of interest for future studies. The theoretical investigation of the 
considered variational formulations, as well as problems with nonlocal 
boundary conditions in general, is far from being a trivial task. Usually, 
the application of the Lax--Milgram theorem in the standard $H^1$ space is 
not possible. In some cases, the existence and uniqueness of the solution in 
the standard function spaces can be proved using embedding theorems. 
Otherwise, the introduction of special (weighted) function spaces is 
necessary. The reader is referred to~\cite{Avalishvili2011,Avalishvili2018} 
for some examples of such investigations, while the theoretical investigation 
of the variational formulations considered in this paper is left as a future 
work. 

\section*{CRediT authorship contribution statement}
\addcontentsline{toc}{section}{CRediT authorship contribution statement}

\textbf{Svaj\={u}nas Sajavičius:} Conceptualization, Methodology, Software, Investigation, Writing - original draft. \textbf{Thomas Takacs:} Validation, Formal analysis, Writing - review \& editing.

\section*{Acknowledgements}
\addcontentsline{toc}{section}{Acknowledgements}
Thomas Takacs was partially supported by the \href{http://dx.doi.org/10.13039/501100002428}{Austrian Science Fund (FWF)} and 
the government of Upper Austria through the project P 30926-NBL. This support 
is gratefully acknowledged.

\bibliographystyle{model1-num-names}
\bibliography{main}

\begin{thebibliography}{62}
\expandafter\ifx\csname natexlab\endcsname\relax\def\natexlab#1{#1}\fi
\providecommand{\url}[1]{\texttt{#1}}
\providecommand{\href}[2]{#2}
\providecommand{\path}[1]{#1}
\providecommand{\DOIprefix}{doi:}
\providecommand{\ArXivprefix}{arXiv:}
\providecommand{\URLprefix}{URL: }
\providecommand{\Pubmedprefix}{pmid:}
\providecommand{\doi}[1]{\href{https://doi.org/#1}{\path{#1}}}
\providecommand{\Pubmed}[1]{\href{pmid:#1}{\path{#1}}}
\providecommand{\bibinfo}[2]{#2}
\ifx\xfnm\relax \def\xfnm[#1]{\unskip,\space#1}\fi
%Type = Article
\bibitem[{Belytschko et~al.(1994)Belytschko, Lu, and Gu}]{Belytschko1994}
\bibinfo{author}{T.~Belytschko}, \bibinfo{author}{Y.~Y. Lu},
  \bibinfo{author}{L.~Gu},
\newblock \bibinfo{title}{Element-free {G}alerkin methods},
\newblock \bibinfo{journal}{Int. J. Numer. Methods Engrg.} \bibinfo{volume}{37}
  (\bibinfo{year}{1994}) \bibinfo{pages}{229--256}.
%Type = Article
\bibitem[{Liu et~al.(1995)Liu, Jun, and Zhang}]{Liu1995}
\bibinfo{author}{W.~K. Liu}, \bibinfo{author}{S.~Jun}, \bibinfo{author}{Y.~F.
  Zhang},
\newblock \bibinfo{title}{Reproducing kernel particle methods},
\newblock \bibinfo{journal}{Internat. J. Numer. Methods Fluids}
  \bibinfo{volume}{20} (\bibinfo{year}{1995}) \bibinfo{pages}{1081--1106}.
%Type = Article
\bibitem[{Bonet and Lok(1999)}]{Bonet1999}
\bibinfo{author}{J.~Bonet}, \bibinfo{author}{T.-S.~L. Lok},
\newblock \bibinfo{title}{Variational and momentum preservation aspects of
  smooth particle hydrodynamic formulations},
\newblock \bibinfo{journal}{Comput. Methods Appl. Mech. Engrg.}
  \bibinfo{volume}{180} (\bibinfo{year}{1999}) \bibinfo{pages}{97--115}.
%Type = Article
\bibitem[{Dumont et~al.(2006)Dumont, Goubet, Ha-Duong, and Villon}]{Dumont2006}
\bibinfo{author}{S.~Dumont}, \bibinfo{author}{O.~Goubet},
  \bibinfo{author}{T.~Ha-Duong}, \bibinfo{author}{P.~Villon},
\newblock \bibinfo{title}{Meshfree methods and boundary conditions},
\newblock \bibinfo{journal}{Int. J. Numer. Methods Engrg.} \bibinfo{volume}{67}
  (\bibinfo{year}{2006}) \bibinfo{pages}{989--1011}.
%Type = Article
\bibitem[{Zhu and Atluri(1998)}]{Zhu1998}
\bibinfo{author}{T.~Zhu}, \bibinfo{author}{S.~N. Atluri},
\newblock \bibinfo{title}{A modified collocation method and a penalty
  formulation for enforcing the essential boundary conditions in the element
  free {G}alerkin method},
\newblock \bibinfo{journal}{Comput. Mech.} \bibinfo{volume}{21}
  (\bibinfo{year}{1998}) \bibinfo{pages}{211--222}.
%Type = Article
\bibitem[{Bonet and Kulasegaram(2000)}]{Bonet2000}
\bibinfo{author}{J.~Bonet}, \bibinfo{author}{S.~Kulasegaram},
\newblock \bibinfo{title}{Correction and stabilization of smooth particle
  hydrodynamics methods with applications in metal forming simulations},
\newblock \bibinfo{journal}{Int. J. Numer. Methods Engrg.} \bibinfo{volume}{47}
  (\bibinfo{year}{2000}) \bibinfo{pages}{1189--1214}.
%Type = Article
\bibitem[{Cho et~al.(2008)Cho, Song, and Choi}]{Cho2008}
\bibinfo{author}{J.~Y. Cho}, \bibinfo{author}{Y.~M. Song},
  \bibinfo{author}{Y.~H. Choi},
\newblock \bibinfo{title}{Boundary locking induced by penalty enforcement of
  essential boundary conditions in mesh-free methods},
\newblock \bibinfo{journal}{Comput. Methods Appl. Mech. Engrg.}
  \bibinfo{volume}{197} (\bibinfo{year}{2008}) \bibinfo{pages}{1167--1183}.
%Type = Article
\bibitem[{Chu and Moran(1995)}]{Chu1995}
\bibinfo{author}{Y.~A. Chu}, \bibinfo{author}{B.~Moran},
\newblock \bibinfo{title}{A computational model for nucleation of solid--solid
  phase transformations},
\newblock \bibinfo{journal}{Model. Simul. Mater. Sci. Engrg.}
  \bibinfo{volume}{3} (\bibinfo{year}{1995}) \bibinfo{pages}{455--471}.
%Type = Article
\bibitem[{Gosz and Liu(1996)}]{Gosz1996}
\bibinfo{author}{J.~Gosz}, \bibinfo{author}{W.~K. Liu},
\newblock \bibinfo{title}{Admissible approximations for essential boundary
  conditions in the reproducing kernel particle method},
\newblock \bibinfo{journal}{Comput. Mech.} \bibinfo{volume}{19}
  (\bibinfo{year}{1996}) \bibinfo{pages}{120--135}.
%Type = Article
\bibitem[{G\"unther and Liu(1998)}]{Guenther1998}
\bibinfo{author}{F.~C. G\"unther}, \bibinfo{author}{W.~K. Liu},
\newblock \bibinfo{title}{Implementation of boundary conditions for meshless
  methods},
\newblock \bibinfo{journal}{Comput. Methods Appl. Mech. Engrg.}
  \bibinfo{volume}{163} (\bibinfo{year}{1998}) \bibinfo{pages}{205--230}.
%Type = Article
\bibitem[{Chen and Wang(2000)}]{Chen2000}
\bibinfo{author}{J.-S. Chen}, \bibinfo{author}{H.-P. Wang},
\newblock \bibinfo{title}{New boundary condition treatments in meshfree
  computation of contact problems},
\newblock \bibinfo{journal}{Comput. Methods Appl. Mech. Engrg.}
  \bibinfo{volume}{187} (\bibinfo{year}{2000}) \bibinfo{pages}{441--468}.
%Type = Article
\bibitem[{Wagner and Liu(2000)}]{Wagner2000}
\bibinfo{author}{G.~J. Wagner}, \bibinfo{author}{W.~K. Liu},
\newblock \bibinfo{title}{Application of essential boundary conditions in
  mesh-free methods: a corrected collocation method},
\newblock \bibinfo{journal}{Int. J. Numer. Methods Engrg.} \bibinfo{volume}{47}
  (\bibinfo{year}{2000}) \bibinfo{pages}{1367--1379}.
%Type = Article
\bibitem[{Sukumar(2004)}]{Sukumar2004}
\bibinfo{author}{N.~Sukumar},
\newblock \bibinfo{title}{Construction of polygonal interpolants: a maximum
  entropy approach},
\newblock \bibinfo{journal}{Int. J. Numer. Methods Engrg.} \bibinfo{volume}{61}
  (\bibinfo{year}{2004}) \bibinfo{pages}{2159--2181}.
%Type = Article
\bibitem[{Oh and Jeong(2009)}]{Oh2009}
\bibinfo{author}{H.-S. Oh}, \bibinfo{author}{J.~W. Jeong},
\newblock \bibinfo{title}{Almost everywhere partition of unity to deal with
  essential boundary conditions in meshless methods},
\newblock \bibinfo{journal}{Comput. Methods Appl. Mech. Engrg.}
  \bibinfo{volume}{198} (\bibinfo{year}{2009}) \bibinfo{pages}{3299--3312}.
%Type = Article
\bibitem[{Belytschko et~al.(1995)Belytschko, Organ, and
  Krongauz}]{Belytschko1995}
\bibinfo{author}{T.~Belytschko}, \bibinfo{author}{D.~Organ},
  \bibinfo{author}{Y.~Krongauz},
\newblock \bibinfo{title}{A coupled finite element--element-free {G}alerkin
  method},
\newblock \bibinfo{journal}{Comput. Mech.} \bibinfo{volume}{17}
  (\bibinfo{year}{1995}) \bibinfo{pages}{186--195}.
%Type = Article
\bibitem[{Krongauz and Belytschko(1996)}]{Krongauz1996}
\bibinfo{author}{Y.~Krongauz}, \bibinfo{author}{T.~Belytschko},
\newblock \bibinfo{title}{Enforcement of essential boundary conditions in
  meshless approximations using finite elements},
\newblock \bibinfo{journal}{Comput. Methods Appl. Mech. Engrg.}
  \bibinfo{volume}{131} (\bibinfo{year}{1996}) \bibinfo{pages}{133--145}.
%Type = Article
\bibitem[{Huerta and Fern{\'a}ndez-M{\'e}ndez(2000)}]{Huerta2000}
\bibinfo{author}{A.~Huerta}, \bibinfo{author}{S.~Fern{\'a}ndez-M{\'e}ndez},
\newblock \bibinfo{title}{Enrichment and coupling of the finite element and
  meshless methods},
\newblock \bibinfo{journal}{Int. J. Numer. Methods Engrg.} \bibinfo{volume}{48}
  (\bibinfo{year}{2000}) \bibinfo{pages}{1615--1636}.
%Type = Article
\bibitem[{Wagner and Liu(2001)}]{Wagner2001}
\bibinfo{author}{G.~J. Wagner}, \bibinfo{author}{W.~K. Liu},
\newblock \bibinfo{title}{Hierarchical enrichment for bridging scales and
  mesh-free boundary conditions},
\newblock \bibinfo{journal}{Int. J. Numer. Methods Engrg.} \bibinfo{volume}{50}
  (\bibinfo{year}{2001}) \bibinfo{pages}{507--524}.
%Type = Article
\bibitem[{Fougeron and Aubry(2019)}]{Fougeron2019}
\bibinfo{author}{G.~Fougeron}, \bibinfo{author}{D.~Aubry},
\newblock \bibinfo{title}{Imposition of boundary conditions for elliptic
  equations in the context of non boundary fitted meshless methods},
\newblock \bibinfo{journal}{Comput. Methods Appl. Mech. Engrg.}
  \bibinfo{volume}{343} (\bibinfo{year}{2019}) \bibinfo{pages}{506--529}.
%Type = Article
\bibitem[{Hughes et~al.(2005)Hughes, Cottrell, and Bazilevs}]{Hughes2005}
\bibinfo{author}{T.~J.~R. Hughes}, \bibinfo{author}{J.~A. Cottrell},
  \bibinfo{author}{Y.~Bazilevs},
\newblock \bibinfo{title}{Isogeometric analysis: {CAD}, finite elements,
  {NURBS}, exact geometry and mesh refinement},
\newblock \bibinfo{journal}{Comput. Methods Appl. Mech. Engrg.}
  \bibinfo{volume}{194} (\bibinfo{year}{2005}) \bibinfo{pages}{4135--4195}.
%Type = Book
\bibitem[{Cottrell et~al.(2009)Cottrell, Hughes, and Bazilevs}]{Cottrell2009}
\bibinfo{author}{J.~A. Cottrell}, \bibinfo{author}{T.~J.~R. Hughes},
  \bibinfo{author}{Y.~Bazilevs}, \bibinfo{title}{Isogeometric Analysis: Toward
  Integration of {CAD} and {FEA}}, \bibinfo{publisher}{Wiley},
  \bibinfo{year}{2009}.
%Type = Article
\bibitem[{Beir\~{a}o~da Veiga et~al.(2014)Beir\~{a}o~da Veiga, Buffa, Sangalli,
  and V\'{a}zquez}]{BeiraodaVeiga2014}
\bibinfo{author}{L.~Beir\~{a}o~da Veiga}, \bibinfo{author}{A.~Buffa},
  \bibinfo{author}{G.~Sangalli}, \bibinfo{author}{R.~V\'{a}zquez},
\newblock \bibinfo{title}{Mathematical analysis of variational isogeometric
  methods},
\newblock \bibinfo{journal}{Acta Numer.} \bibinfo{volume}{23}
  (\bibinfo{year}{2014}) \bibinfo{pages}{157--287}.
%Type = Inproceedings
\bibitem[{Hughes and Sangalli(2017)}]{Hughes2017}
\bibinfo{author}{T.~J.~R. Hughes}, \bibinfo{author}{G.~Sangalli},
\newblock \bibinfo{title}{Mathematics of isogeometric analysis: A conspectus},
\newblock in: \bibinfo{editor}{E.~Stein}, \bibinfo{editor}{R.~de~Borst},
  \bibinfo{editor}{T.~J.~R. Hughes} (Eds.), \bibinfo{booktitle}{Encyclopedia of
  Computational Mechanics, Second Edition, Volume 1, Part 2. Fundamentals},
  \bibinfo{year}{2017}.
%Type = Inproceedings
\bibitem[{Hughes et~al.(2018)Hughes, Sangalli, and Tani}]{Hughes2018}
\bibinfo{author}{T.~J.~R. Hughes}, \bibinfo{author}{G.~Sangalli},
  \bibinfo{author}{M.~Tani},
\newblock \bibinfo{title}{Isogeometric analysis: {M}athematical and
  implementational aspects, with applications},
\newblock in: \bibinfo{booktitle}{Splines and {PDE}s: {F}rom Approximation
  Theory to Numerical Linear Algebra}, volume \bibinfo{volume}{2219} of
  \textit{\bibinfo{series}{Lecture Notes in Mathematics}},
  \bibinfo{publisher}{Springer, Cham}, \bibinfo{year}{2018}, pp.
  \bibinfo{pages}{237--315}.
%Type = Article
\bibitem[{Nguyen et~al.(2015)Nguyen, Anitescu, Bordas, and
  Rabczuk}]{Nguyen2015}
\bibinfo{author}{V.~P. Nguyen}, \bibinfo{author}{C.~Anitescu},
  \bibinfo{author}{S.~P.~A. Bordas}, \bibinfo{author}{T.~Rabczuk},
\newblock \bibinfo{title}{Isogeometric analysis: {A}n overview and computer
  implementation aspects},
\newblock \bibinfo{journal}{Math. Comput. Simulation} \bibinfo{volume}{117}
  (\bibinfo{year}{2015}) \bibinfo{pages}{89--116}.
%Type = Article
\bibitem[{Costantini et~al.(2010)Costantini, Manni, Pelosi, and
  Sampoli}]{Costantini2010}
\bibinfo{author}{P.~Costantini}, \bibinfo{author}{C.~Manni},
  \bibinfo{author}{F.~Pelosi}, \bibinfo{author}{M.~L. Sampoli},
\newblock \bibinfo{title}{Quasi-interpolation in isogeometric analysis based on
  generalized {B}-splines},
\newblock \bibinfo{journal}{Comput. Aided Geom. Design} \bibinfo{volume}{27}
  (\bibinfo{year}{2010}) \bibinfo{pages}{656--668}.
%Type = Article
\bibitem[{Wang and Xuan(2010)}]{Wang2010b}
\bibinfo{author}{D.~Wang}, \bibinfo{author}{J.~Xuan},
\newblock \bibinfo{title}{An improved {NURBS}-based isogeometric analysis with
  enhanced treatment of essential boundary conditions},
\newblock \bibinfo{journal}{Comput. Methods Appl. Mech. Engrg.}
  \bibinfo{volume}{199} (\bibinfo{year}{2010}) \bibinfo{pages}{2425--2436}.
%Type = Article
\bibitem[{Chen et~al.(2011)Chen, Mo, and Gong}]{Chen2011b}
\bibinfo{author}{T.~Chen}, \bibinfo{author}{R.~Mo}, \bibinfo{author}{Z.~W.
  Gong},
\newblock \bibinfo{title}{Imposing essential boundary conditions in
  isogeometric analysis with {N}itsche's method},
\newblock \bibinfo{journal}{Applied Mechanics and Materials}
  \bibinfo{volume}{121--126} (\bibinfo{year}{2011})
  \bibinfo{pages}{2779--2783}.
%Type = Inproceedings
\bibitem[{Mitchell et~al.(2011)Mitchell, Govindjee, and Taylor}]{Mitchell2011}
\bibinfo{author}{T.~J. Mitchell}, \bibinfo{author}{S.~Govindjee},
  \bibinfo{author}{R.~L. Taylor},
\newblock \bibinfo{title}{A method for enforcement of {D}irichlet boundary
  conditions in isogeometric analysis},
\newblock in: \bibinfo{editor}{D.~Mueller-Hoeppe},
  \bibinfo{editor}{S.~L\"{o}ehnert}, \bibinfo{editor}{S.~Reese} (Eds.),
  \bibinfo{booktitle}{Recent Developments and Innovative Applications in
  Computational Mechanics}, \bibinfo{publisher}{Springer Berlin Heidelberg},
  \bibinfo{year}{2011}, pp. \bibinfo{pages}{283--293}.
%Type = Article
\bibitem[{Govindjee et~al.(2012)Govindjee, Strain, Mitchell, and
  Taylor}]{Govindjee2012}
\bibinfo{author}{S.~Govindjee}, \bibinfo{author}{J.~Strain},
  \bibinfo{author}{T.~J. Mitchell}, \bibinfo{author}{R.~L. Taylor},
\newblock \bibinfo{title}{Convergence of an efficient local least-squares
  fitting method for bases with compact support},
\newblock \bibinfo{journal}{Comput. Methods Appl. Mech. Engrg.}
  \bibinfo{volume}{213--216} (\bibinfo{year}{2012}) \bibinfo{pages}{84--92}.
%Type = Article
\bibitem[{Bazilevs and Hughes(2007)}]{Bazilevs2007}
\bibinfo{author}{Y.~Bazilevs}, \bibinfo{author}{T.~J.~R. Hughes},
\newblock \bibinfo{title}{Weak imposition of {D}irichlet boundary conditions in
  fluid mechanics},
\newblock \bibinfo{journal}{Comput. Fluids} \bibinfo{volume}{36}
  (\bibinfo{year}{2007}) \bibinfo{pages}{12--26}.
%Type = Article
\bibitem[{Bazilevs et~al.(2007)Bazilevs, Michler, Calo, and
  Hughes}]{Bazilevs2007b}
\bibinfo{author}{Y.~Bazilevs}, \bibinfo{author}{C.~Michler},
  \bibinfo{author}{V.~M. Calo}, \bibinfo{author}{T.~J.~R. Hughes},
\newblock \bibinfo{title}{Weak {D}irichlet boundary conditions for wall-bounded
  turbulent flows},
\newblock \bibinfo{journal}{Comput. Methods Appl. Mech. Engrg.}
  \bibinfo{volume}{196} (\bibinfo{year}{2007}) \bibinfo{pages}{4853--4862}.
%Type = Article
\bibitem[{Ruess et~al.(2014)Ruess, Schillinger, \"Ozcan, and Rank}]{Ruess2014}
\bibinfo{author}{M.~Ruess}, \bibinfo{author}{D.~Schillinger},
  \bibinfo{author}{A.~I. \"Ozcan}, \bibinfo{author}{E.~Rank},
\newblock \bibinfo{title}{Weak coupling for isogeometric analysis of
  non-matching and trimmed multi-patch geometries},
\newblock \bibinfo{journal}{Comput. Methods Appl. Mech. Engrg.}
  \bibinfo{volume}{269} (\bibinfo{year}{2014}) \bibinfo{pages}{46--71}.
%Type = Article
\bibitem[{Schillinger et~al.(2012)Schillinger, Ded\`{e}, Scott, Evans, Borden,
  Rank, and Hughes}]{Schillinger2012}
\bibinfo{author}{D.~Schillinger}, \bibinfo{author}{L.~Ded\`{e}},
  \bibinfo{author}{M.~A. Scott}, \bibinfo{author}{J.~A. Evans},
  \bibinfo{author}{M.~J. Borden}, \bibinfo{author}{E.~Rank},
  \bibinfo{author}{T.~J.~R. Hughes},
\newblock \bibinfo{title}{An isogeometric design-through-analysis methodology
  based on adaptive hierarchical refinement of {NURBS}, immersed boundary
  methods, and {T}-spline {CAD} surfaces},
\newblock \bibinfo{journal}{Comput. Methods Appl. Mech. Engrg.}
  \bibinfo{volume}{249--252} (\bibinfo{year}{2012}) \bibinfo{pages}{116--150}.
%Type = Article
\bibitem[{Rank et~al.(2012)Rank, Ruess, Kollmannsberger, Schillinger, and
  D\"uster}]{Rank2012}
\bibinfo{author}{E.~Rank}, \bibinfo{author}{M.~Ruess},
  \bibinfo{author}{S.~Kollmannsberger}, \bibinfo{author}{D.~Schillinger},
  \bibinfo{author}{A.~D\"uster},
\newblock \bibinfo{title}{Geometric modeling, isogeometric analysis and the
  finite cell method},
\newblock \bibinfo{journal}{Comput. Methods Appl. Mech. Engrg.}
  \bibinfo{volume}{249--252} (\bibinfo{year}{2012}) \bibinfo{pages}{104--115}.
%Type = Article
\bibitem[{Nitti et~al.(2020)Nitti, Kiendl, Reali, and {de~Tullio}}]{Nitti2020}
\bibinfo{author}{A.~Nitti}, \bibinfo{author}{J.~Kiendl},
  \bibinfo{author}{A.~Reali}, \bibinfo{author}{M.~D. {de~Tullio}},
\newblock \bibinfo{title}{An immersed-boundary/isogeometric method for
  fluid--structure interaction involving thin shells},
\newblock \bibinfo{journal}{Comput. Methods Appl. Mech. Engrg}
  \bibinfo{volume}{364} (\bibinfo{year}{2020}) \bibinfo{pages}{112977}.
%Type = Article
\bibitem[{Fern\'{a}ndez-M\'{e}ndez and Huerta(2004)}]{Fernandez-Mendez2004}
\bibinfo{author}{S.~Fern\'{a}ndez-M\'{e}ndez}, \bibinfo{author}{A.~Huerta},
\newblock \bibinfo{title}{Imposing essential boundary conditions in mesh-free
  methods},
\newblock \bibinfo{journal}{Comput. Methods Appl. Mech. Engrg.}
  \bibinfo{volume}{193} (\bibinfo{year}{2004}) \bibinfo{pages}{1257--1275}.
%Type = Inproceedings
\bibitem[{Huerta et~al.(2017)Huerta, Belytschko, Fern\'{a}ndez-M\'{e}ndez,
  Rabczuk, Zhuang, and Arroyo}]{Huerta2017}
\bibinfo{author}{A.~Huerta}, \bibinfo{author}{T.~Belytschko},
  \bibinfo{author}{S.~Fern\'{a}ndez-M\'{e}ndez}, \bibinfo{author}{T.~Rabczuk},
  \bibinfo{author}{X.~Zhuang}, \bibinfo{author}{M.~Arroyo},
\newblock \bibinfo{title}{Meshfree methods},
\newblock in: \bibinfo{editor}{E.~Stein}, \bibinfo{editor}{R.~de~Borst},
  \bibinfo{editor}{T.~J.~R. Hughes} (Eds.), \bibinfo{booktitle}{Encyclopedia of
  Computational Mechanics, Second Edition, Part 2. Fundamentals},
  \bibinfo{year}{2017}, pp. \bibinfo{pages}{1--38}.
%Type = Article
\bibitem[{Embar et~al.(2010)Embar, Dolbow, and Harari}]{Embar2010}
\bibinfo{author}{A.~Embar}, \bibinfo{author}{J.~Dolbow},
  \bibinfo{author}{I.~Harari},
\newblock \bibinfo{title}{Imposing {D}irichlet boundary conditions with
  {N}itsche's method and spline-based finite elements},
\newblock \bibinfo{journal}{Int. J. Numer. Methods Engrg.} \bibinfo{volume}{83}
  (\bibinfo{year}{2010}) \bibinfo{pages}{877--898}.
%Type = Article
\bibitem[{Day(1985)}]{Day1985b}
\bibinfo{author}{W.~A. Day},
\newblock \bibinfo{title}{Existence of a property of solutions of the heat
  equation subject to linear thermoelasticity and other theories},
\newblock \bibinfo{journal}{Quart. Appl. Math.} \bibinfo{volume}{40}
  (\bibinfo{year}{1985}) \bibinfo{pages}{319--330}.
%Type = Article
\bibitem[{Day(1992)}]{Day1992}
\bibinfo{author}{W.~A. Day},
\newblock \bibinfo{title}{Parabolic equations and thermodynamics},
\newblock \bibinfo{journal}{Quart. Appl. Math.} \bibinfo{volume}{50}
  (\bibinfo{year}{1992}) \bibinfo{pages}{523--533}.
%Type = Article
\bibitem[{Hazanee and Lesnic(2014)}]{Hazanee2014}
\bibinfo{author}{A.~Hazanee}, \bibinfo{author}{D.~Lesnic},
\newblock \bibinfo{title}{Determination of a time-dependent coefficient in the
  bioheat equation},
\newblock \bibinfo{journal}{Int. J. Mech. Sci.} \bibinfo{volume}{88}
  (\bibinfo{year}{2014}) \bibinfo{pages}{259--266}.
%Type = Article
\bibitem[{D\'{\i}az et~al.(1998)D\'{\i}az, Padial, and Rakotoson}]{Diaz1998}
\bibinfo{author}{J.~I. D\'{\i}az}, \bibinfo{author}{J.~F. Padial},
  \bibinfo{author}{J.~M. Rakotoson},
\newblock \bibinfo{title}{Mathematical treatment of the magnetic confinement in
  a current carrying stellarator},
\newblock \bibinfo{journal}{Nonlinear Anal.} \bibinfo{volume}{34}
  (\bibinfo{year}{1998}) \bibinfo{pages}{857--887}.
%Type = Article
\bibitem[{Glotov et~al.(2016)Glotov, Hames, Meirc, and Ngoma}]{Glotov2016}
\bibinfo{author}{D.~Glotov}, \bibinfo{author}{W.~E. Hames},
  \bibinfo{author}{A.~J. Meirc}, \bibinfo{author}{S.~Ngoma},
\newblock \bibinfo{title}{An integral constrained parabolic problem with
  applications in thermochronology},
\newblock \bibinfo{journal}{Comp. Math. Appl.} \bibinfo{volume}{71}
  (\bibinfo{year}{2016}) \bibinfo{pages}{2301--2312}.
%Type = Article
\bibitem[{Madenci et~al.(2018)Madenci, Dorduncu, Barut, and Phan}]{Madenci2018}
\bibinfo{author}{E.~Madenci}, \bibinfo{author}{M.~Dorduncu},
  \bibinfo{author}{A.~Barut}, \bibinfo{author}{N.~Phan},
\newblock \bibinfo{title}{Weak form of peridynamics for nonlocal essential and
  natural boundary conditions},
\newblock \bibinfo{journal}{Comput. Methods Appl. Mech. Engrg.}
  \bibinfo{volume}{337} (\bibinfo{year}{2018}) \bibinfo{pages}{598--631}.
%Type = Article
\bibitem[{Ainsworth(2001)}]{Ainsworth2001}
\bibinfo{author}{M.~Ainsworth},
\newblock \bibinfo{title}{Essential boundary conditions and multi-point
  constraints in finite element analysis},
\newblock \bibinfo{journal}{Comput. Methods Appl. Mech. Engrg.}
  \bibinfo{volume}{190} (\bibinfo{year}{2001}) \bibinfo{pages}{6323--6339}.
%Type = Article
\bibitem[{Jendele and \v{C}ervenka(2009)}]{Jendele2009}
\bibinfo{author}{L.~Jendele}, \bibinfo{author}{J.~\v{C}ervenka},
\newblock \bibinfo{title}{On the solution of multi-point constraints --
  {A}pplication to {FE} analysis of reinforced concrete structures},
\newblock \bibinfo{journal}{Comput. Struct.} \bibinfo{volume}{87}
  (\bibinfo{year}{2009}) \bibinfo{pages}{970--980}.
%Type = Article
\bibitem[{Kashiwabara et~al.(2015)Kashiwabara, Colciago, Ded\`e, and
  Quarteroni}]{Kashiwabara2015}
\bibinfo{author}{T.~Kashiwabara}, \bibinfo{author}{C.~M. Colciago},
  \bibinfo{author}{L.~Ded\`e}, \bibinfo{author}{A.~Quarteroni},
\newblock \bibinfo{title}{Well-posedness, regularity, and convergence analysis
  of the finite element approximation of a generalized {R}obin boundary value
  problem},
\newblock \bibinfo{journal}{SIAM J. Numer. Anal.} \bibinfo{volume}{53}
  (\bibinfo{year}{2015}) \bibinfo{pages}{105--126}.
%Type = Misc
\bibitem[{Bercovier and Soloveichik(2015)}]{Bercovier2015}
\bibinfo{author}{M.~Bercovier}, \bibinfo{author}{I.~Soloveichik},
  \bibinfo{title}{Overlapping non matching meshes domain decomposition method
  in isogeometric analysis}, \bibinfo{year}{2015}.
  \href{http://arxiv.org/abs/1502.03756}{\tt arXiv:1502.03756}.
%Type = Article
\bibitem[{Kargaran et~al.(2019)Kargaran, J\"uttler, Kleiss, Mantzaflaris, and
  Takacs}]{Kargaran2019}
\bibinfo{author}{S.~Kargaran}, \bibinfo{author}{B.~J\"uttler},
  \bibinfo{author}{S.~K. Kleiss}, \bibinfo{author}{A.~Mantzaflaris},
  \bibinfo{author}{T.~Takacs},
\newblock \bibinfo{title}{Overlapping multi-patch structures in isogeometric
  analysis},
\newblock \bibinfo{journal}{Comput. Methods Appl. Mech. Engrg.}
  \bibinfo{volume}{356} (\bibinfo{year}{2019}) \bibinfo{pages}{325--353}.
%Type = Article
\bibitem[{Avalishvili et~al.(2011)Avalishvili, Avalishvili, and
  Gordeziani}]{Avalishvili2011}
\bibinfo{author}{G.~Avalishvili}, \bibinfo{author}{M.~Avalishvili},
  \bibinfo{author}{D.~Gordeziani},
\newblock \bibinfo{title}{On a nonlocal problem with integral boundary
  conditions for a multidimensional elliptic equation},
\newblock \bibinfo{journal}{Appl. Math. Lett.} \bibinfo{volume}{24}
  (\bibinfo{year}{2011}) \bibinfo{pages}{566--571}.
%Type = Article
\bibitem[{Avalishvili et~al.(2018)Avalishvili, Avalishvili, and
  Miara}]{Avalishvili2018}
\bibinfo{author}{G.~Avalishvili}, \bibinfo{author}{M.~Avalishvili},
  \bibinfo{author}{B.~Miara},
\newblock \bibinfo{title}{Nonclassical problem with integral boundary
  conditions for elliptic system},
\newblock \bibinfo{journal}{Complex Var. Elliptic Equ.} \bibinfo{volume}{63}
  (\bibinfo{year}{2018}) \bibinfo{pages}{836--853}.
%Type = Article
\bibitem[{Sajavi\v{c}ius(2014)}]{Sajavicius2014a}
\bibinfo{author}{S.~Sajavi\v{c}ius},
\newblock \bibinfo{title}{Radial basis function method for a multidimensional
  linear elliptic equation with nonlocal boundary conditions},
\newblock \bibinfo{journal}{Comput. Math. Appl.} \bibinfo{volume}{67}
  (\bibinfo{year}{2014}) \bibinfo{pages}{1407--1420}.
%Type = Article
\bibitem[{Sajavi\v{c}ius(2016)}]{Sajavicius2016}
\bibinfo{author}{S.~Sajavi\v{c}ius},
\newblock \bibinfo{title}{Radial basis function collocation method for an
  elliptic problem with nonlocal multipoint boundary condition},
\newblock \bibinfo{journal}{Engrg. Anal. Bound. Elem.} \bibinfo{volume}{67}
  (\bibinfo{year}{2016}) \bibinfo{pages}{164--172}.
%Type = Book
\bibitem[{Piegl and Tiller(1997)}]{Piegl1997}
\bibinfo{author}{L.~Piegl}, \bibinfo{author}{W.~Tiller}, \bibinfo{title}{The
  {NURBS} Book}, Monographs in Visual Communications, \bibinfo{edition}{2nd}
  ed., \bibinfo{publisher}{Springer Berlin Heidelberg}, \bibinfo{year}{1997}.
%Type = Book
\bibitem[{Rogers(2001)}]{Rogers2001}
\bibinfo{author}{D.~F. Rogers}, \bibinfo{title}{An Introduction to {NURBS}:
  With Historical Perspective}, \bibinfo{publisher}{Morgan Kaufmann},
  \bibinfo{year}{2001}.
%Type = Book
\bibitem[{de~{B}oor(2001)}]{deBoor2001}
\bibinfo{author}{C.~de~{B}oor}, \bibinfo{title}{A Practical Guide to Splines},
  \bibinfo{edition}{revised} ed., \bibinfo{publisher}{Springer-Verlag New
  York}, \bibinfo{year}{2001}.
%Type = Article
\bibitem[{Auricchio et~al.(2010)Auricchio, Beir\~{a}o~da Veiga, Hughes, Reali,
  and Sangalli}]{Auricchio2010}
\bibinfo{author}{F.~Auricchio}, \bibinfo{author}{L.~Beir\~{a}o~da Veiga},
  \bibinfo{author}{T.~J.~R. Hughes}, \bibinfo{author}{A.~Reali},
  \bibinfo{author}{G.~Sangalli},
\newblock \bibinfo{title}{Isogeometric collocation methods},
\newblock \bibinfo{journal}{Math. Models Methods Appl. Sci.}
  \bibinfo{volume}{20} (\bibinfo{year}{2010}) \bibinfo{pages}{2075}.
%Type = Inproceedings
\bibitem[{Reali and Hughes(2015)}]{Reali2015b}
\bibinfo{author}{A.~Reali}, \bibinfo{author}{T.~J.~R. Hughes},
\newblock \bibinfo{title}{An introduction to isogeometric collocation methods},
\newblock in: \bibinfo{editor}{G.~Beer}, \bibinfo{editor}{S.~Bordas} (Eds.),
  \bibinfo{booktitle}{Isogeometric Methods for Numerical Simulation}, volume
  \bibinfo{volume}{561} of \textit{\bibinfo{series}{CISM International Centre
  for Mechanical Sciences}}, \bibinfo{publisher}{Springer Vienna},
  \bibinfo{year}{2015}, pp. \bibinfo{pages}{173--204}.
%Type = Article
\bibitem[{Marussig and Hughes(2018)}]{Marussig2018}
\bibinfo{author}{B.~Marussig}, \bibinfo{author}{T.~J.~R. Hughes},
\newblock \bibinfo{title}{A review of trimming in isogeometric analysis:
  {C}hallenges, data exchange and simulation aspects},
\newblock \bibinfo{journal}{Arch. Computat. Methods. Engrg.}
  \bibinfo{volume}{25} (\bibinfo{year}{2018}) \bibinfo{pages}{1059--1127}.
%Type = Article
\bibitem[{V\'{a}zquez(2016)}]{Vazquez2016}
\bibinfo{author}{R.~V\'{a}zquez},
\newblock \bibinfo{title}{A new design for the implementation of isogeometric
  analysis in {O}ctave and {M}atlab: {GeoPDEs} 3.0},
\newblock \bibinfo{journal}{Comput. Math. Appl.} \bibinfo{volume}{72}
  (\bibinfo{year}{2016}) \bibinfo{pages}{523--554}.
%Type = Misc
\bibitem[{Geo(2016)}]{GeoPDEs2016}
\bibinfo{title}{{GeoPDEs}: a package for {I}sogeometric {A}nalysis in {M}atlab
  and {O}ctave}, \bibinfo{year}{2016}. \URLprefix
  \url{http://rafavzqz.github.io/geopdes/}.

\end{thebibliography}
\addcontentsline{toc}{section}{References}

\end{document}